\documentclass[10pt]{amsart}
\usepackage{amscd}
\usepackage{amssymb}
\usepackage[all]{xy}
\usepackage{graphicx}
\author{Amnon Yekutieli}
\title{Continuous and Twisted $\mrm{L}_{\infty}$ Morphisms}
\date{6 September 2005}
\address{Department of Mathematics,
Ben Gurion University, Be'er Sheva 84105, Israel}
\email{amyekut@math.bgu.ac.il}
\thanks{{\em Mathematics Subject Classification} 2000.
Primary: 53D55; Secondary: 13D10, 13N10, 13J10.}
\keywords{$\mrm{L}_{\infty}$ morphisms, deformation quantization.}

\newtheorem{thm}[equation]{Theorem}
\newtheorem{cor}[equation]{Corollary}
\newtheorem{prop}[equation]{Proposition}
\newtheorem{lem}[equation]{Lemma}
\theoremstyle{definition}
\newtheorem{dfn}[equation]{Definition}
\newtheorem{rem}[equation]{Remark}
\newtheorem{exa}[equation]{Example}

\numberwithin{equation}{section}

\newcommand{\iso}{\xrightarrow{\simeq}}
\newcommand{\inj}{\hookrightarrow}

\newcommand{\xar}{\xrightarrow}
\newcommand{\opn}{\operatorname}
\newcommand{\cat}[1]{\operatorname{\mathsf{#1}}}

\newcommand{\ol}{\overline}

\newcommand{\rmitem}[1]{\item[\text{\textup{(#1)}}]}
\newcommand{\mfrak}[1]{\mathfrak{#1}}
\newcommand{\mcal}[1]{\mathcal{#1}}

\newcommand{\mrm}[1]{\mathrm{#1}}
\newcommand{\mbb}[1]{\mathbb{#1}}

\newcommand{\smfrac}[2]{{\textstyle \frac{#1}{#2}}}
\newcommand{\tup}[1]{\textup{#1}}
\newcommand{\bsym}[1]{\boldsymbol{#1}}
\newcommand{\boplus}{\bigoplus\nolimits}

\newcommand{\til}[1]{\tilde{#1}}
\newcommand{\what}[1]{\widehat{#1}}

\newcommand{\set}[1]{\{ #1 \}}

\newcommand{\K}{\mbb{K} \hspace{0.05em}}
\renewcommand{\d}{\mrm{d}}
\newcommand{\bwedge}{{\textstyle \bigwedge}}
\newcommand{\sprod}{{\textstyle \prod}}
\newcommand{\hatotimes}[1]{\, \what{\otimes}_{#1} \,}

\begin{document}

\begin{abstract}
The purpose of this paper is to develop a suitable notion of 
continuous $\mrm{L}_{\infty}$ morphism between DG Lie algebras, 
and to study twists of such morphisms. 
\end{abstract}

\maketitle

\setcounter{section}{-1}
\section{Introduction}

Let $\K$ be a field containing $\mbb{R}$. 
Consider two DG Lie algebras associated to the polynomial algebra
$\K[\bsym{t}] := \K[t_1, \ldots, t_n]$. The first is the algebra of 
{\em poly derivations} $\mcal{T}_{\mrm{poly}}(\K[\bsym{t}])$,
and the second is the algebra of {\em poly differential operators} 
$\mcal{D}_{\mrm{poly}}(\K[\bsym{t}])$.
A very important result of M. Kontsevich \cite{Ko1}, 
known as the Formality 
Theorem, gives an explicit formula for an $\mrm{L}_{\infty}$ 
quasi-isomorphism 
\[ \mcal{U} : \mcal{T}_{\mrm{poly}}(\K[\bsym{t}])
\to \mcal{D}_{\mrm{poly}}(\K[\bsym{t}]) . \]

Here is the main result of our paper. 

\begin{thm} \label{thm0.1}
Assume $\mbb{R} \subset \K$.
Let $A = \boplus_{i \geq 0} A^i$ be a super-commutative 
associative unital complete DG algebra in 
$\cat{Dir} \cat{Inv} \cat{Mod} \K$. Consider the 
induced continuous $A$-multilinear $\mrm{L}_{\infty}$ morphism
\[ \mcal{U}_A : 
A \hatotimes{} \mcal{T}_{\mrm{poly}}(\K[[\bsym{t}]]) \to
A \hatotimes{} \mcal{D}_{\mrm{poly}}(\K[[\bsym{t}]]) . \]
Suppose 
$\omega \in 
A^1 \hatotimes{} \mcal{T}^0_{\mrm{poly}}(\K[[\bsym{t}]])$
is a solution of the Maurer-Cartan equation in 
$A \hatotimes{} \mcal{T}_{\mrm{poly}}(\K[[\bsym{t}]])$.
Define 
$\omega' := (\partial^1 \mcal{U}_A)(\omega) \in
A^1 \hatotimes{} \mcal{D}^0_{\mrm{poly}}(\K[[\bsym{t}]])$.
Then $\omega'$ is a solution of the Maurer-Cartan equation in 
$A \hatotimes{} \mcal{D}_{\mrm{poly}}(\K[[\bsym{t}]])$,
and there is continuous $A$-multilinear $\mrm{L}_{\infty}$ 
quasi-isomorphism 
\[ \mcal{U}_{A, \omega} :
\big( A \hatotimes{} \mcal{T}_{\mrm{poly}}(\K[[\bsym{t}]]) 
\big)_{\omega} \to
\big( A \hatotimes{} \mcal{D}_{\mrm{poly}}(\K[[\bsym{t}]]) 
\big)_{\omega'} \]
whose Taylor coefficients are
\[ (\partial^j \mcal{U}_{A, \omega})(\alpha) :=
\sum_{k \geq 0} \smfrac{1}{(j+k)!} 
(\partial^{j+k} \mcal{U}_{A})(\omega^k \wedge \alpha)  \]
for 
$\alpha \in \sprod^j
\bigl( A \hatotimes{} \mcal{T}_{\mrm{poly}}(\K[[\bsym{t}]])
\bigr)$.
\end{thm}

Below is an outline of the paper, in which we mention the 
various terms appearing in the theorem. 

In Section 1 we 
develop the theory of dir-inv modules. A dir-inv structure on a 
$\K$-module $M$ is a generalization of an adic topology. The 
category of dir-inv modules and continuous homomorphisms is denoted 
by $\cat{Dir} \cat{Inv} \cat{Mod} \K$. 
The concept of dir-inv module, 
and related complete tensor product $\hatotimes{}$, are quite 
flexible, and are particulary well-suited for infinitely generated 
modules. Among other things we introduce the notion of DG Lie 
algebra in $\cat{Dir} \cat{Inv} \cat{Mod} \K$.

Section 2 concentrates on poly differential operators. The results 
here are mostly generalizations of material from \cite{EGA-IV}.

In Section 3 we review the coalgebra approach to 
$\mrm{L}_{\infty}$ morphisms. The notions of continuous, 
$A$-multilinear and twisted 
$\mrm{L}_{\infty}$ morphism are defined. The main 
result of this section is Theorem \ref{thm9.2}.

In Section 4 we recall the Kontsevich Formality Theorem. 
By combining it with Theorem \ref{thm9.2} we deduce Theorem 
\ref{thm0.1} (repeated as Theorem \ref{thm4.10}). In Theorem 
\ref{thm0.1} the DG Lie algebras
$A \hatotimes{} \mcal{T}_{\mrm{poly}}(\K[[\bsym{t}]])$
and 
$A \hatotimes{} \mcal{D}_{\mrm{poly}}(\K[[\bsym{t}]])$
are the $A$-multilinear extensions of
$\mcal{T}_{\mrm{poly}}(\K[[\bsym{t}]])$ and
$\mcal{D}_{\mrm{poly}}(\K[[\bsym{t}]])$ respectively,
and
$\big( A \hatotimes{} \mcal{T}_{\mrm{poly}}(\K[[\bsym{t}]]) 
\big)_{\omega}$
and
$\big( A \hatotimes{} \mcal{D}_{\mrm{poly}}(\K[[\bsym{t}]]) 
\big)_{\omega'}$
are their twists. The $\mrm{L}_{\infty}$ morphism
 $\mcal{U}_{A}$ is the continuous
$A$-multilinear extension of $\mcal{U}$, and 
$\mcal{U}_{A, \omega}$ is its twist.

Theorem \ref{thm0.1} is used in \cite{Ye2}, in which we study 
deformation quantization of algebraic varieties. 

\medskip \noindent
\textbf{Acknowledgments.} The author wishes to thank Vladimir 
Hinich, Bernhard Keller and James Stasheff for their assistance.

\section{Dir-Inv Modules}
\label{sec1} 

We begin the paper with a generalization of the notion of adic 
topology. In this section $\K$ is a commutative base ring, and $C$ 
is a commutative $\K$-algebra. 
The category $\cat{Mod} C$ is abelian and has 
direct and inverse limits. Unless specified otherwise, all 
limits are taken in $\cat{Mod} C$.

\begin{dfn}
\begin{enumerate}
\item Let $M \in \cat{Mod} C$. An {\em inv module 
structure} on $M$ is an inverse system 
$\{ \mrm{F}^i M \}_{i \in \mbb{N}}$ of $C$-submodules
of $M$. The pair 
$(M, \{ \mrm{F}^i M \}_{i \in \mbb{N}})$
is called an {\em inv $C$-module}.
\item Let $(M, \{ \mrm{F}^i M \}_{i \in \mbb{N}})$ and
$(N, \{ \mrm{F}^i N \}_{i \in \mbb{N}})$ be two inv $C$-modules.
A function $\phi : M \to N$ ($C$-linear or not) 
is said to be {\em continuous}
if for every $i \in \mbb{N}$ there exists $i' \in \mbb{N}$ such 
that $\phi(\mrm{F}^{i'} M) \subset \mrm{F}^i N$.
\item Define $\cat{Inv} \cat{Mod} C$
to be the category whose objects are the inv $C$-modules, 
and whose morphisms are the continuous $C$-linear homomorphisms.  
\end{enumerate}
\end{dfn}

We do not assume that the canonical homomorphism 
$M \to \lim_{\leftarrow i} M / \mrm{F}^i M$
is surjective nor injective. There is a full embedding 
$\cat{Mod} C \inj \cat{Inv} \cat{Mod} C$,
$M \mapsto (M, \{ \ldots, 0, 0 \})$. 
If $(M, \{ \mrm{F}^i M \}_{i \in \mbb{N}})$ and
$(N, \{ \mrm{F}^i N \}_{i \in \mbb{N}})$ 
are two inv $C$-modules then 
$M \oplus N$ is an inv module, with inverse system of 
submodules 
$\mrm{F}^i (M \oplus N) := \mrm{F}^i M \oplus \mrm{F}^i N$.
Thus $\cat{Inv} \cat{Mod} C$ is a $C$-linear additive 
category.

Let $(M, \{ \mrm{F}^i M \}_{i \in \mbb{N}})$ be an inv $C$-module,
let $M', M''$ be two $C$-modules, and suppose 
$\phi : M' \to M$ and $\psi : M \to M''$ are $C$-linear 
homomorphisms. We get induced inv module structures on $M'$ 
and $M''$ by defining 
$\mrm{F}^i M' := \phi^{-1}(\mrm{F}^i M)$ and 
$\mrm{F}^i M'' := \psi(\mrm{F}^i M)$.

Recall that a {\em directed set} is a partially ordered set $J$ 
with the property that for any $j_1, j_2 \in J$ there exists 
$j_3 \in J$ such that $j_1, j_2 \leq j_3$.

\begin{dfn}
\begin{enumerate}
\item Let $M \in \cat{Mod} C$. A {\em 
dir-inv module structure} on $M$ is a direct system 
$\{ \mrm{F}_j M \}_{j \in J}$ of $C$-submodules of $M$,
indexed by a nonempty directed set $J$, 
together with an inv module structure on 
each $\mrm{F}_j M$, such that for every $j_1 \leq j_2$ the 
inclusion $\mrm{F}_{j_1} M \inj \mrm{F}_{j_2} M$
is continuous. The pair $(M, \{ \mrm{F}_j M \}_{j \in J})$
is called a {\em dir-inv $C$-module}.
\item Let $(M, \{ \mrm{F}_j M \})_{j \in J}$ and
$(N, \{ \mrm{F}_k N \}_{k \in K})$ 
be two dir-inv $C$-modules.
A function $\phi : M \to N$ ($C$-linear or not) 
 is said to be {\em continuous}
if for every $j \in J$ there exists $k \in K$ 
such that $\phi(\mrm{F}_j M) \subset \mrm{F}_k N$, and
$\phi : \mrm{F}_j M \to \mrm{F}_k N$ is a continuous 
function between these two inv $C$-modules.
\item Define $\cat{Dir} \cat{Inv} \cat{Mod} C$
to be the category whose objects are the dir-inv $C$-modules, 
and whose morphisms are the continuous $C$-linear 
homomorphisms.  
\end{enumerate}
\end{dfn}

There is no requirement that the canonical homomorphism
$\lim_{j \to} \mrm{F}_j M \to M$ 
will be surjective. An inv $C$-module $M$ is endowed with the
dir-inv module structure $\{ \mrm{F}_j M \}_{j \in J}$, where 
$J := \{ 0 \}$ and $\mrm{F}_0 M := M$.
Thus we get a full embedding 
$\cat{Inv} \cat{Mod} C \inj 
\cat{Dir} \cat{Inv} \cat{Mod} C$. 
Given two dir-inv $C$-modules 
$(M, \{ \mrm{F}_j M \})_{j \in J}$ and
$(N, \{ \mrm{F}_k N \}_{k \in K})$, 
we make $M \oplus N$ into a dir-inv module as 
follows. The directed set is $J \times K$, with the component-wise 
partial order, and the direct system of inv modules is 
$\mrm{F}_{(j, k)} (M \oplus N) :=
\mrm{F}_j M \oplus \mrm{F}_k N$.
The condition $J \neq \emptyset$ in part (1) of the definition 
ensures that the zero module $0 \in \cat{Mod} C$ is an initial 
object in $\cat{Dir} \cat{Inv} \cat{Mod} C$.
So $\cat{Dir} \cat{Inv} \cat{Mod} C$ is a $C$-linear 
additive category.

Let $(M, \{ \mrm{F}_j M \}_{j \in J})$ be a dir-inv $C$-module,
let $M', M''$ be two $C$-modules, and suppose 
$\phi : M' \to M$ and $\psi : M \to M''$ are $C$-linear 
homomorphisms. We get induced dir-inv module structures 
$\{ \mrm{F}_j M' \}_{j \in J}$ and
$\{ \mrm{F}_j M'' \}_{j \in J}$
on $M'$ and $M''$ as follows. Define
$\mrm{F}_j(M') := \phi^{-1}(\mrm{F}_j M)$ and
$\mrm{F}_j M'' := \psi(\mrm{F}_j M)$,
which have induced inv module structures via the homomorphisms
$\phi : \mrm{F}_j M' \to \mrm{F}_j M$ and
$\psi : \mrm{F}_j M \to \mrm{F}_j M''$.

\begin{dfn} \label{dfn1.3}
\begin{enumerate}
\item An inv $C$-module 
$(M, \{ \mrm{F}^i M \}_{i \in \mbb{N}})$ is called 
{\em discrete} if \linebreak $\mrm{F}^i M = 0$ for $i \gg 0$.
\item An inv $C$-module 
$(M, \{ \mrm{F}^i M \}_{i \in \mbb{N}})$
is called {\em complete} if the canonical homomorphism
$M \to  \lim_{\leftarrow i} M / \mrm{F}^i M$ is bijective.
\item A dir-inv $C$-module $M$ is called 
{\em complete} (resp.\ {\em discrete}) if it isomorphic, in 
$\cat{Dir} \cat{Inv} \cat{Mod} C$,
to a dir-inv module 
$(N, \{ \mrm{F}_j N \}_{j \in J})$,
where all the inv modules $\mrm{F}_j N$ are complete 
(resp.\ discrete) as defined above, and the canonical homomorphism
$\lim_{j \to} \mrm{F}_j N \to N$ is bijective.
\item A dir-inv $C$-module $M$ is called 
{\em trivial} if it isomorphic, in 
$\cat{Dir} \cat{Inv} \cat{Mod} C$,
to an object of $\cat{Mod} C$,
via the embedding 
$\cat{Mod} C \inj \cat{Dir} \cat{Inv} \cat{Mod} C$.
\end{enumerate}
\end{dfn}

Note that $M$ is a trivial dir-inv module iff it is 
isomorphic, in $\cat{Dir} \cat{Inv} \cat{Mod} C$, 
to a discrete inv module. There do exist discrete dir-inv modules 
that are not trivial dir-inv modules; see Example \ref{exa2.4}.
It is easy to see that if $M$ is a discrete dir-inv module 
then it is also complete. 

The base ring $\K$ is endowed with the inv structure 
$\{ \ldots, 0, 0 \}$, so it is a trivial dir-inv $\K$-module. 
But the $\K$-algebra $C$ could have more interesting dir-inv 
structures (cf.\ Example \ref{exa2.2}).

If $f^* : C \to C'$ is a homomorphism of $\K$-algebras, then there 
is a functor
$f_* : \cat{Dir} \cat{Inv} \cat{Mod} C' \to 
\cat{Dir} \cat{Inv} \cat{Mod} C$. 
In particular any dir-inv $C$-module is a dir-inv $\K$-module.

\begin{dfn} \label{dfn2.5}
\begin{enumerate}
\item Given an inv $C$-module 
$(M, \{ \mrm{F}^i M \}_{i \in \mbb{N}})$ its completion 
is the inv $C$-module 
$(\what{M}, \{ \mrm{F}^i \what{M} \}_{i \in \mbb{N}})$,
defined as follows: 
$\what{M} := \lim_{\leftarrow i} M / \mrm{F}^i M$
and
$\mrm{F}^i \what{M} := \opn{Ker}(\what{M} \to 
M / \mrm{F}^i M)$.
Thus we obtain an additive endofunctor $M \mapsto \what{M}$ of 
$\cat{Inv} \cat{Mod} C$.
\item Given a dir-inv $C$-module 
$(M, \{ \mrm{F}_j M \}_{j \in J})$ its completion
is the dir-inv $C$-module 
$(\what{M}, \{ \mrm{F}_j \what{M} \}_{j \in J})$ 
defined as follows. For any $j \in J$ let
$\what{\mrm{F}_j M}$ be the completion of the inv $C$-module 
$\mrm{F}_j M$, as defined above. Then let 
$\what{M} := \lim_{j \to} \what{\mrm{F}_j M}$
and 
$\mrm{F}_j \what{M} := \opn{Im}(\what{\mrm{F}_j M} \to \what{M})$. 
Thus we obtain an additive endofunctor $M \mapsto \what{M}$ of 
$\cat{Dir} \cat{Inv} \cat{Mod} C$.
\end{enumerate}
\end{dfn}

An inv $C$-module $M$ is complete iff the functorial 
homomorphism 
$M \to \what{M}$ is an isomorphism; and of course
$\what{M}$ is complete. For a dir-inv $C$-module
$M$ there is in general no functorial homomorphism between 
$M$ and $\what{M}$, and we do not know if 
$\what{M}$ is complete. Nonetheless:

\begin{prop} \label{prop1.7}
Suppose $M \in \cat{Dir} \cat{Inv} \cat{Mod} C$
is complete. Then there is an isomorphism 
$M \cong \what{M}$ in 
$\cat{Dir} \cat{Inv} \cat{Mod} C$.
This isomorphism is functorial.
\end{prop}

\begin{proof}
For any dir-inv module 
$(M, \{ \mrm{F}_j M \}_{j \in J})$ 
let's define
$M' := \lim_{j \to} \mrm{F}_j M$.
So $(M', \{ \mrm{F}_j M \}_{j \in J})$
is a dir-inv module, and there are functorial morphisms 
$M' \to M$ and
$M' \to \what{M}$. 
If $M$ is complete then both these morphisms are 
isomorphisms.
\end{proof}

Suppose $\{ M_k \}_{k \in K}$ is a collection of dir-inv 
modules, indexed by a set $K$. There is an induced
dir-inv module structure 
on $M := \boplus_{k \in K}  M_k$, constructed as follows. For any 
$k$ let us denote by $\{ \mrm{F}_j M_{k} \}_{j \in J_k}$
the dir-inv structure of $M_{k}$; so that each
$\mrm{F}_j M_{k}$ is an inv module. 
For each finite subset $K_0 \subset K$ let
$J_{K_0} := \prod_{k \in K_0} J_k$,
made into a directed set by component-wise partial order. 
Define $J := \coprod_{K_0} J_{K_0}$, where $K_0$ runs over the 
finite subsets of $K$.
For two finite subsets $K_0 \subset K_1$, and two elements
$\bsym{j}_0 = \{ j_{0,k} \}_{k \in K_0} \in J_{K_0}$
and
$\bsym{j}_1 = \{ j_{1,k} \}_{k \in K_1} \in J_{K_1}$
we declare that $\bsym{j}_0 \leq \bsym{j}_1$ if
$j_{0,k} \leq j_{1,k}$ for all $k \in K_0$.
This makes $J$ into a directed set. Now for any 
$\bsym{j} = \{ j_{k} \}_{k \in K_0} \in J_{K_0} \subset J$
let $\mrm{F}_{\bsym{j}} M := \boplus_{k \in K_0} 
\mrm{F}_{j_k} M_k$, which is an inv module. The dir-inv structure 
on $M$ is $\{ F_{\bsym{j}} M \}_{\bsym{j} \in J}$.

\begin{prop} \label{prop1.15}
Let $\{ M_k \}_{k \in K}$ be a collection of dir-inv 
$C$-modules, and let $M := \boplus_{k \in K}  M_k$, endowed with 
the induced dir-inv structure.
\begin{enumerate}
\item $M$ is a coproduct of $\{ M_k \}_{k \in K}$ 
in the category $\cat{Dir} \cat{Inv} \cat{Mod} C$.
\item There is a functorial isomorphism
$\what{M} \cong \boplus_{k \in K} \what{M}_k$.
\end{enumerate}
\end{prop}

\begin{proof}
(1) is obvious. For (2) we note that both $\what{M}$ and
$\boplus_{k \in K} \what{M}_k$ are direct limits for the 
direct system 
$\{ \what{M}_{\bsym{j}} \}_{\bsym{j} \in J}$.
\end{proof}

Suppose $\set{M_k}_{k \in \mbb{N}}$ is a collection of inv 
$C$-modules. For each $k$ let 
$\set{ \mrm{F}^i M_k }_{i \in \mbb{N}}$ be the inv structure of 
$M_k$. Then 
$M := \prod_{k \in \mbb{N}} M_k$ is an inv module, with inv 
structure 
$\mrm{F}^i M := \big( \prod_{k > i} M_k \big) \times 
\big( \prod_{k \leq i} \mrm{F}^i M_k \big)$.
Next let $\set{M_k}_{k \in \mbb{N}}$ be a collection of dir-inv 
$C$-modules, and for each $k$ let 
$\set{ \mrm{F}_j M_k }_{j \in J_k}$ be the dir-inv structure of 
$M_k$. Then there is an induced dir-inv structure on 
$M := \prod_{k \in \mbb{N}} M_k$.
Define a directed set $J :=  \prod_{k \in \mbb{N}} J_k$,
with component-wise partial order. For any 
$\bsym{j} = \set{j_k}_{k \in \mbb{N}} \in J$
define 
$\mrm{F}_{\bsym{j}} M := \prod_{k \in \mbb{N}}
\mrm{F}_{j_k} M_k$, 
which is an inv $C$-module as explained above. The dir-inv 
structure on $M$ is 
$\set{ \mrm{F}_{\bsym{j}} M }_{\bsym{j} \in J}$.

\begin{prop}
Let $\set{M_k}_{k \in \mbb{N}}$ be a collection of dir-inv 
$C$-modules, and let $M := \prod_{k \in \mbb{N}} M_k$, endowed 
with the induced dir-inv structure. 
Then $M$ is a product of $\set{M_k}_{k \in \mbb{N}}$ in 
$\cat{Dir} \cat{Inv} \cat{Mod} C$.
\end{prop}

\begin{proof}
All we need to consider is continuity. 
First assume that all the $M_k$ are inv $C$-modules.
Let's denote by $\pi_k : M \to M_k$ the projection. 
For each $k, i \in \mbb{N}$ and $i' \geq \opn{max}(i, k)$ we have
$\pi_{k}(\mrm{F}^{i'} M) = \mrm{F}^i M_k$.
This shows that the $\pi_k$ are continuous. 
Suppose $L$ is an inv $C$-module and $\phi_k : L \to M_k$ 
are morphisms in $\cat{Inv} \cat{Mod} C$. 
For any $i \in \mbb{N}$
there exists $i' \in \mbb{N}$ such that 
$\phi_k(\mrm{F}^{i'} L) \subset \mrm{F}^i M_k$
for all $k \leq i$. Therefore the homomorphism $\phi : L \to M$
with components $\phi_k$ is continuous.

Now let $M_k$ be dir-inv $C$-modules, with dir-inv structures
$\set{ \mrm{F}_j M_k }_{j \in J_k}$.
For any $\bsym{j} = \set{j_k} \in J$ one has
$\pi_k(\mrm{F}_{\bsym{j}} M) = \mrm{F}_{j_k} M_k$,
and as shown above
$\pi_k :  \mrm{F}_{\bsym{j}} M \to \mrm{F}_{j_k} M_k$
is continuous. Given a dir-inv module $L$ and 
and morphisms $\phi_k : L \to M_k$ in 
$\cat{Dir} \cat{Inv} \cat{Mod} C$, we have to prove that 
$\phi : L \to M$ is continuous. 
Let $\set{ \mrm{F}_j L }_{j \in J_L}$ be the dir-inv structure of 
$L$. Take any $j \in J_L$. Since $\phi_k$ is continuous, 
there exists some $j_k \in J_k$ such that 
$\phi_k(\mrm{F}_j L) \subset \mrm{F}_{j_k} M_k$. But then
$\phi(\mrm{F}_j L) \subset \mrm{F}_{\bsym{j}} M$
for $\bsym{j} := \set{j_k}_{k \in \mbb{N}}$, and 
by the previous paragraph 
$\phi : \mrm{F}_j L \to \mrm{F}_{\bsym{j}} M$
is continuous.
\end{proof}

The following examples should help to clarify the notion of dir-inv 
module.

\begin{exa} \label{exa2.2}
Let $\mfrak{c}$ be an ideal in $C$.
Then each finitely generated $C$-module $M$ has an inv structure 
$\set{ \mrm{F}^i M }_{i \in \mbb{N}}$,
where we define the submodules
$\mrm{F}^i M := \mfrak{c}^{i+1} M$. This is called the 
{\em $\mfrak{c}$-adic inv structure}.
Any $C$-module $M$ has a dir-inv structure 
$\{ \mrm{F}_j M \}_{j \in J}$, which is  
the collection of finitely generated $C$-submodules of $M$, 
directed by inclusion, and each $\mrm{F}_j M$ is given the 
$\mfrak{c}$-adic inv structure. We get a fully faithful functor
$\cat{Mod} C \to \cat{Dir} \cat{Inv} \cat{Mod} C$.
This dir-inv module structure on $M$ is called the 
{\em $\mfrak{c}$-adic dir-inv structure}.

In case $C$ is noetherian and $\mfrak{c}$-adically complete,
then the finitely generated modules are complete as 
inv $C$-modules, and hence all modules are complete as dir-inv 
modules.
\end{exa}

\begin{exa} \label{exa1.5}  
Suppose $(M, \{ \mrm{F}^i M \}_{i \in \mbb{N}})$ is an inv 
$C$-module, and $\{ i_k \}_{k \in \mbb{N}}$ is a nondecreasing 
sequence in $\mbb{N}$ with 
$\opn{lim}_{k \to \infty} i_k = \infty$. 
Then $\{ \mrm{F}^{i_k} M \}_{k \in \mbb{N}}$
is a new inv structure on $M$, yet the identity map
$(M, \{ \mrm{F}^i M \}_{i \in \mbb{N}}) \to
(M, \{ \mrm{F}^{i_k} M \}_{k \in \mbb{N}})$
is an isomorphism in $\cat{Inv} \cat{Mod} C$.

A similar modification can be done for dir-inv modules. Suppose 
$(M, \{ \mrm{F}_j M \}_{j \in J})$ is a dir-inv $C$-module, and
$J' \subset J$ is a subset that is cofinal in 
$J$. Then $\{ \mrm{F}_j M \}_{j \in J'}$ is a new dir-inv 
structure on $M$, yet the identity map
$(M, \{ \mrm{F}_j M \}_{j \in J}) \to$ \linebreak
$(M, \{ \mrm{F}_j M \}_{j \in J'})$
is an isomorphism in $\cat{Dir} \cat{Inv} \cat{Mod} C$.
\end{exa}

\begin{exa} \label{exa2.4}
Let $M$ be the free $\K$-module 
with basis $\{ e_p \}_{p \in \mbb{N}}$; so 
$M =$ \linebreak $\boplus_{p \in \mbb{N}} \K e_p$ in $\cat{Mod} \K$.
We put on $M$ the inv module structure 
$\{ \mrm{F}^i M \}_{i \in \mbb{N}}$ with $\mrm{F}^i M := 0$
for all $i$. 
Let $N$ be the same $\K$-module as $M$, but put on it the inv 
module structure $\{ \mrm{F}^i N \}_{i \in \mbb{N}}$ with
$\mrm{F}^i N := \boplus_{p=i}^{\infty} \K e_p$. 
Also let $L$ be the $\K$-module $M$, but put on it the 
dir-inv module structure $\{ \mrm{F}_j L \}_{j \in \mbb{N}}$, with
$\mrm{F}_j L := \boplus_{p=0}^j \K e_p$ the discrete inv module
whose inv structure is $\{ \ldots, 0, 0 \}$.
Both $L$ and $M$ are discrete and complete as dir-inv 
$\K$-modules, and 
$\what{N} \cong \prod_{p \in \mbb{N}} \K e_p$.
The dir-inv module $M$ is trivial. 
$L$ is not a trivial dir-inv $\K$-module, because it is not 
isomorphic in $\cat{Dir} \cat{Inv} \cat{Mod} \K$
to any inv module. 
The identity maps $L \to M \to N$ are continuous. The only 
continuous $\K$-linear homomorphisms $M \to L$ are those with 
finitely generated images. 
\end{exa}

\begin{rem}
In the situation of the previous example, 
suppose we put on the three modules $L, M, N$ genuine 
$\K$-linear topologies, using the limiting processes and starting 
from the discrete topology. Namely $M$, $N / \mrm{F}^i N$ and 
$\mrm{F}_j L$ get the discrete topologies; 
$L \cong \lim_{j \to} \mrm{F}_j L$
gets the $\lim_{\to}$ topology; and 
$N \subset \lim_{\leftarrow i} N / \mrm{F}^i N$ gets the 
$\lim_{\leftarrow}$ topology (as in \cite[Section 1.1]{Ye1}).
Then $L$ and $M$ become the same discrete topological module, and 
$\what{N}$ is the topological completion of $N$.
We see that the notion of a dir-inv structure is more 
subtle than that of a topology, even though a similar 
language is used.
\end{rem}

\begin{exa}
Suppose $\K$ is a field, and let $M := \K$, the free module of 
rank $1$. Up to isomorphism in 
$\cat{Dir} \cat{Inv} \cat{Mod} \K$, $M$ has three distinct dir-inv 
module structures. We can denote them by $M_1, M_2, M_3$ in such a 
way that the identity maps $M_1 \to M_2 \to M_3$ 
are continuous. The only continuous $\K$-linear homomorphisms
$M_i \to M_j$ with $i > j$ are the zero homomorphisms. $M_2$ is 
the trivial dir-inv structure, and it is the only interesting one 
(the others are ``pathological'').
\end{exa}

\begin{exa} \label{exa1.4}
Suppose $M = \boplus_{p \in \mbb{Z}} M^p$ is a graded $C$-module. 
The grading induces a dir-inv structure on $M$, with 
$J := \mbb{N}$,
$\mrm{F}_j M := \boplus_{p = -j}^{\infty} M^p$, and
$\mrm{F}^i \mrm{F}_j M := \boplus_{p = -j + i}^{\infty} M^p$.
The completion satisfies 
$\what{M} \cong \big( \prod_{p \geq 0} M^p \big) \oplus
\big( \boplus_{p < 0} M^p \big)$ in
$\cat{Dir} \cat{Inv} \cat{Mod} C$, 
where each $M^p$ has the trivial dir-inv module structure.
\end{exa}
 
It makes sense to talk about convergence of sequences in a dir-inv 
module. Suppose $(M, \{ \mrm{F}^i M \}_{i \in \mbb{N}})$ 
is an inv $C$-module and  
$\{ m_i \}_{i \in \mbb{N}}$ is a sequence in $M$.
We say that $\lim_{i \to \infty} m_i = 0$
if for every $i_0$ there is some $i_1$ such that 
$\{ m_i \}_{i \geq i_1} \subset \mrm{F}_{i_0} M$. If
$(M, \{ \mrm{F}_j M \}_{j \in J})$ is a dir-inv module
and $\{ m_i \}_{i \in \mbb{N}}$ is a sequence in 
$M$, then we say that $\lim_{i \to \infty} m_i = 0$
if there exist some $j$ and $i_1$ such that 
$\{ m_i \}_{i \geq i_1} \subset \mrm{F}_j M$,
and $\lim_{i \to \infty} m_i = 0$ in the inv module 
$\mrm{F}_j M$. Having defined 
$\lim_{i \to \infty} m_i = 0$, it is clear how to define
$\lim_{i \to \infty} m_i = m$
and
$\sum_{i = 0}^{\infty} m_i = m$. 
Also the notion of Cauchy sequence is clear. 

\begin{prop}
Assume $M$ is a complete dir-inv $C$-module. Then any Cauchy 
sequence in $M$ has a unique limit.
\end{prop}

\begin{proof}
Consider a Cauchy sequence $\{ m_i \}_{i \in \mbb{N}}$ in $M$.
Convergence is an invariant of isomorphisms in 
$\cat{Dir} \cat{Inv} \cat{Mod} C$.
By Definition \ref{dfn1.3} we may assume that in the 
dir-inv structure $\{ \mrm{F}_j M \}_{j \in J}$ 
of $M$ each inv module $\mrm{F}_j M$ is complete. By passing to 
the sequence  $\{ m_i - m_{i_1} \}_{i \in \mbb{N}}$
for suitable $i_1$, 
we can also assume the sequence is contained in one of the inv 
modules $\mrm{F}_j M$. Thus we reduce to the case of convergence 
in a complete inv module, which is standard. 
\end{proof}

Let $(M, \{ \mrm{F}^i M \}_{i \in \mbb{N}})$ and
$(N, \{ \mrm{F}^i N \}_{i \in \mbb{N}})$
be two inv $C$-modules. We make 
$M \otimes_{C} N$ into an inv module by 
defining 
\[ \mrm{F}^i (M \otimes_{C} N) :=
\opn{Im} \big( (M \otimes_{C} \mrm{F}^i N)
\oplus (\mrm{F}^i M \otimes_{C} N) \to
M \otimes_{C} N \big) . \]
For two dir-inv $C$-modules
$(M, \{ \mrm{F}_j M \}_{j \in J})$ and
$(N, \{ \mrm{F}_k N \}_{k \in  K})$,
we put on $M \otimes_{C} N$ 
the dir-inv module structure 
$\{ \mrm{F}_{(j,k)} (M \otimes_{C} N) \}
_{(j,k) \in J \times K}$,
where 
\[ \mrm{F}_{(j,k)} (M \otimes_{C} N) :=
\opn{Im}(\mrm{F}_j M \otimes_{C} \mrm{F}_k N \to
M \otimes_{C} N) . \]

\begin{dfn} \label{dfn2.7}
Given $M, N \in \cat{Dir} \cat{Inv} \cat{Mod} C$ we define 
$N \hatotimes{C} M$
to be the completion of the dir-inv $C$-module
$N \otimes_{C} M$.
\end{dfn}

\begin{exa}
Let's examine the behavior of the dir-inv modules $L, M, N$ from 
Example \ref{exa2.4} with respect to complete tensor 
product. There is an isomorphism 
$L \otimes_{\K} N \cong \boplus_{p \in \mbb{N}} N$
in $\cat{Dir} \cat{Inv} \cat{Mod} \K$,
so according to Proposition \ref{prop1.15}(2) there is also an
isomorphism
$L \hatotimes{\K} N \cong \boplus_{p \in \mbb{N}} \what{N}$
in $\cat{Dir} \cat{Inv} \cat{Mod} \K$.
On the other hand $M \otimes_{\K} N$ is an inv $\K$-module
with inv structure 
$\mrm{F}^i (M \otimes_{\K} N) = M \otimes_{\K} \mrm{F}^i N$,
so $M \hatotimes{\K} N \cong \prod_{p \in \mbb{N}} M$
in $\cat{Dir} \cat{Inv} \cat{Mod} \K$.
The series $\sum_{p=0}^{\infty} e_p \otimes e_p$
converges in $M \hatotimes{\K} N$, but not in 
$L \hatotimes{\K} N$.
\end{exa}

A {\em graded object in $\cat{Dir} \cat{Inv} \cat{Mod} C$},
or a {\em graded dir-inv $C$-module},
is an object $M \in \cat{Dir} \cat{Inv} \cat{Mod} C$
of the form 
$M = \boplus_{i \in \mbb{Z}} M^i$,
with $M^i \in \cat{Dir} \cat{Inv} \cat{Mod} C$.
According to Proposition \ref{prop1.15} we have
$\what{M} \cong \boplus_{i \in \mbb{Z}} \what{M}^i$.
Given two graded objects
$M = \boplus_{i \in \mbb{Z}} M^i$ and 
$N = \boplus_{i \in \mbb{Z}} N^i$
in $\cat{Dir} \cat{Inv} \cat{Mod} C$, 
the tensor product is also a graded object in 
$\cat{Dir} \cat{Inv} \cat{Mod} C$, with
\[ (M \otimes_{C} N)^i :=
\boplus_{p+q = i} \, M^p \otimes_{C} N^q . \]

In this paper ``algebra'' is taken in the weakest possible sense: 
by $C$-algebra we mean an $C$-module $A$ 
together with an $C$-bilinear function
$\mu_A : A \times A \to A$.
If $A$ is associative, or a Lie algebra, 
then we will specify that. However, ``commutative algebra''
will mean, by default, a commutative associative unital 
$C$-algebra. Another convention is that a homomorphism between 
unital algebras is a unital homomorphism, and a module over a 
unital algebra is a unital module.

\begin{dfn} \label{dfn2.6}
\begin{enumerate}
\item An {\em algebra in} $\cat{Dir} \cat{Inv} \cat{Mod} 
C$ is an object
$A \in$ \linebreak $\cat{Dir} \cat{Inv} \cat{Mod} C$,
together with a continuous $C$-bilinear function
$\mu_A : A \times A \to A$.
\item A {\em differential graded algebra in}
$\cat{Dir} \cat{Inv} \cat{Mod} C$ is a graded object
$A =$ \linebreak $\boplus_{i \in \mbb{Z}} A^i$ 
in $\cat{Dir} \cat{Inv} \cat{Mod} C$,
together with continuous $C$-(bi)linear functions
$\mu_A : A \times A \to A$
and $\mrm{d}_A : A \to A$,
such that $A$ is a differential graded
algebra, in the usual sense, 
with respect to the differential $\mrm{d}_A$ and the 
multiplication $\mu_A$.
\item Let $A$ be an algebra in $\cat{Dir} \cat{Inv} \cat{Mod} C$,
with dir-inv structure $\set{\mrm{F}_j A}_{j \in J}$. 
We say that $A$ is a {\em unital algebra in}
$\cat{Dir} \cat{Inv} \cat{Mod} C$ if it has a unit element
$1_A$ (in the usual sense), such that 
$1_A \in \bigcup_{j \in J} \mrm{F}_j A$.
\end{enumerate}
\end{dfn}

The base ring $\K$, with its trivial dir-inv structure, is a unital 
algebra in \linebreak $\cat{Dir} \cat{Inv} \cat{Mod} \K$.
In item (3) above, the condition 
$1_A \in \bigcup_{j \in J} \mrm{F}_j A$
is equivalent to the ring homomorphism $\K \to A$ being continuous.

We will use the common abbreviation ``DG'' for ``differential 
graded''. An algebra in $\cat{Dir} \cat{Inv} \cat{Mod} C$ can have 
further attributes, such as ``Lie'' or ``associative'', which have 
their usual meanings. 
If $A \in \cat{Inv} \cat{Mod} C$ then we also say it 
is an algebra in $\cat{Inv} \cat{Mod} C$.

\begin{exa} \label{exa1.3}
In the situation of Example \ref{exa2.2},
the $\mfrak{c}$-adic inv structure makes $C$
and $\what{C}$ into unital algebras in $\cat{Inv} \cat{Mod} C$.
\end{exa}

Recall that a graded algebra $A$ is called {\em 
super-commutative} if $b a = (-1)^{i j} a b$ and $c^2 = 0$ for all
$a \in A^i$, $b \in A^j$, $c \in A^k$ and $k$ 
odd. There is no essential difference between left and right DG
$A$-modules.

\begin{prop} \label{prop1.16}  
Let $A$ and $\mfrak{g}$ be DG algebras in 
$\cat{Dir} \cat{Inv} \cat{Mod} C$.
\begin{enumerate}
\item The completion $\what{A}$ is a
DG algebra in $\cat{Dir} \cat{Inv} \cat{Mod} C$.
\item If $A$ is complete, then the canonical isomorphism 
$A \cong \what{A}$ of Proposition \tup{\ref{prop1.7}} 
is an isomorphism of DG algebras.
\item The complete tensor product
$A \what{\otimes}_{C}\, \mfrak{g}$ is a DG algebra in 
$\cat{Dir} \cat{Inv} \cat{Mod} C$.
\item If $A$ is a super-commutative associative unital algebra,
then so is $\what{A}$.
\item If $\mfrak{g}$ is a DG Lie algebra and $A$
is a super-commutative associative unital algebra, then 
$A \what{\otimes}_{C}\, \mfrak{g}$ is a 
DG Lie algebra.
\end{enumerate}
\end{prop} 

\begin{proof}
(1) This is a consequence of a slightly more general fact.
Consider modules
$M_1, \ldots, M_r, N \in \cat{Dir} \cat{Inv} \cat{Mod} C$
and a continuous $C$-multilinear linear function
$\phi : M_1 \times \cdots \times M_r \to N$.
We claim that there is an induced continuous $C$-multilinear 
linear function
$\what{\phi} : \prod_k \what{M}_k \to \what{N}$.
This operation is functorial (w.r.t.\ morphisms $M_k \to M'_k$ and
$N \to N'$), and monoidal (i.e.\ it respects composition in the
$k$th argument with a continuous multilinear function
$\psi : L_{1} \times \cdots \times L_s \to M_k$).

First assume $M_1, \ldots, M_r, N \in \cat{Inv} \cat{Mod} C$,
with inv structures $\set{ \mrm{F}^i M_1 }_{i \in \mbb{N}}$ etc. 
For any $i \in \mbb{N}$ there exists $i' \in \mbb{N}$
such that
$\phi(\prod_k \mrm{F}^{i'} M_k) \subset \mrm{F}^i N$.
Therefore there's an induced continuous $C$-multilinear function 
$\what{\phi} : \prod_k \what{M}_k \to \what{N}$.
It is easy to verify that $\phi \mapsto \what{\phi}$
is functorial and monoidal.

Next consider the general case, i.e.\
$M_1, \ldots, M_r, N \in \cat{Dir} \cat{Inv} \cat{Mod} C$.
Let \linebreak
$\set{ \mrm{F}_j M_k }_{j \in J_{k}}$ be the dir-inv structure of 
$M_k$, and let $\set{ \mrm{F}_j N }_{j \in J_N}$ be the dir-inv 
structure of $N$. By continuity of $\phi$, given 
$(j_1, \ldots, j_r) \in \prod_k J_k$ there exists 
$j' \in J_N$ such that 
$\phi(\prod_k \mrm{F}_{j_k} M_k) \subset \mrm{F}_{j'} N$,
and $\phi : \prod_k \mrm{F}_{j_k} M_k \to \mrm{F}_{j'} N$
is continuous. By the previous paragraph this extends to 
$\what{\phi} : \prod_k \what{\mrm{F}_{j_k} M_k} \to 
\what{\mrm{F}_{j'} N}$.
Passing to the direct limit in 
$(j_1, \ldots, j_r)$ we obtain
$\what{\phi} : \prod_k \what{M}_k \to \what{N}$.
Again this operation is functorial and monoidal.

\medskip \noindent
(2) Let $A' \subset A$ be as in the proof of Proposition 
\ref{prop1.7}. This is a subalgebra. The arguments used in the 
proof of part (1) above show that $A' \to A$ and $A' \to \what{A}$
are algebra homomorphisms.
 
\medskip \noindent
(3) Let us write $\cdot_A$  and $\cdot_{\mfrak{g}}$ for the two 
multiplications, and $\d_A$  and $\d_{\mfrak{g}}$ for the 
differentials. Then $A \otimes_{C} \mfrak{g}$ is a DG algebra 
with multiplication
\[ (a_1 \otimes \gamma_1) \cdot (a_2 \otimes \gamma_2) 
:= (-1)^{i_2 j_1} (a_1 \cdot_A a_2) \otimes 
(\gamma_1 \cdot_{\mfrak{g}} \gamma_2) \]
and differential
\[ \d(a_1 \otimes \gamma_1) := \d_A(a_1) \otimes \gamma_1 +
(-1)^{i_1} a_1 \otimes \d_{\mfrak{g}}(\gamma_1) \]
for $a_k \in A^{i_k}$ and 
$\gamma_k \in \mfrak{g}^{j_k}$. These operations are continuous, 
so $A \otimes_{C} \mfrak{g}$ is a DG algebra in
$\cat{Dir} \cat{Inv} \cat{Mod} C$. Now use part (1).

\medskip \noindent
(4, 5) The various identities (Lie etc.) are 
preserved by $\hatotimes{}$. Definition \ref{dfn2.6}(3) ensures 
that $\what{A}$ has a unit element.
\end{proof}

\begin{dfn} \label{dfn1.4}
Suppose $A$ is a DG super-commutative associative unital
algebra in $\cat{Dir} \cat{Inv} \cat{Mod} C$.
\begin{enumerate}
\item A {\em DG $A$-module in $\cat{Dir} \cat{Inv} \cat{Mod} C$} 
is a graded object $M \in \cat{Dir} \cat{Inv} \cat{Mod} C$, 
together with a continuous $C$-bilinear homomorphism
$A \times M \to M$, which makes $M$ into a DG  
$A$-module in the usual sense.
\item A {\em DG $A$-module Lie algebra in 
$\cat{Dir} \cat{Inv} \cat{Mod} C$} is a DG Lie algebra
$\mfrak{g} \in \cat{Dir} \cat{Inv} \cat{Mod} C$, 
together with a continuous $C$-bilinear homomorphism \linebreak
$A \times \mfrak{g} \to \mfrak{g}$,
such that $\mfrak{g}$ is a DG $A$-module,
and 
\[ [a_1 \gamma_1, a_2 \gamma_2] = (-1)^{i_2 j_1}
a_1 a_2 \,[\gamma_1, \gamma_2] \]
for all $a_k \in A^{i_k}$ and 
$\gamma_k \in \mfrak{g}^{j_k}$.
\end{enumerate}
\end{dfn}

\begin{exa}
If $A$ is a DG super-commutative associative unital 
algebra in $\cat{Dir} \cat{Inv} \cat{Mod} C$, and
$\mfrak{g}$ is a DG Lie algebra in 
$\cat{Dir} \cat{Inv} \cat{Mod} C$, then
$A \hatotimes{C} \mfrak{g}$ is a 
DG $\what{A}$-module Lie algebra in 
$\cat{Dir} \cat{Inv} \cat{Mod} C$.
\end{exa}

Let $A$ be a DG super-commutative associative unital 
algebra in  $\cat{Dir} \cat{Inv} \cat{Mod} C$, and let
$M, N$ be two DG $A$-modules in 
$\cat{Dir} \cat{Inv} \cat{Mod} C$. The tensor product
$M \otimes_{A} N$ is a quotient of $M \otimes_{C} N$,
and as such it has a dir-inv structure. Moreover,
$M \otimes_{A} N$ is a DG $A$-module in 
$\cat{Dir} \cat{Inv} \cat{Mod} C$, and we 
define $M \hatotimes{A} N$ to be its completion, which is 
a DG $\what{A}$-module in $\cat{Dir} \cat{Inv} \cat{Mod} C$.

\begin{prop} \label{prop1.8}
Let $A$ and $B$ be DG super-commutative 
associative unital algebras in $\cat{Dir} \cat{Inv} \cat{Mod} C$, 
and let $A \to B$ be a continuous homomorphism of DG $C$-algebras.
\begin{enumerate}
\item Suppose $M$ is a DG $A$-module in 
$\cat{Dir} \cat{Inv} \cat{Mod} C$. Then 
$B \hatotimes{A} M$ is a DG $\what{B}$-module in 
$\cat{Dir} \cat{Inv} \cat{Mod} C$.
\item Suppose $\mfrak{g}$ is a DG $A$-module Lie algebra in 
$\cat{Dir} \cat{Inv} \cat{Mod} C$. Then 
$B \hatotimes{A} \mfrak{g}$ is a 
DG $\what{B}$-module Lie algebra in 
$\cat{Dir} \cat{Inv} \cat{Mod} C$.
\end{enumerate}
\end{prop}

\begin{proof}
Like Proposition \ref{prop1.16}.
\end{proof}

Suppose $C, C'$ are commutative algebras in 
$\cat{Dir} \cat{Inv} \cat{Mod} \K$, and $f^* : C \to C'$ is a 
continuous $\K$-algebra homomorphism. There are functors
$f^* : \cat{Dir} \cat{Inv} \cat{Mod} C \to
\cat{Dir} \cat{Inv} \cat{Mod} C'$
and
$f^{\what{*}} : \cat{Dir} \cat{Inv} \cat{Mod} C \to
\cat{Dir} \cat{Inv} \cat{Mod} \what{C}'$,
namely
$f^* M := C' \otimes_{C} M$ and 
$f^{\what{*}} M := C' \hatotimes{C} M$. 

Let $M$ and $N$ be two dir-inv $C$-modules. We define
\[ \opn{Hom}^{\mrm{cont}}_{C}(M, N) := 
\opn{Hom}_{\cat{Dir} \cat{Inv} \cat{Mod} C}(M, N) , \]
i.e.\ the $C$-module of continuous $C$-linear homomorphisms.
In general this module has no obvious structure.
However, if $M$ is an inv $C$-module with inv structure
$\{ \mrm{F}^i M \}_{i \in \mbb{N}}$, and $N$ is a discrete inv 
$C$-module, then
\[ \opn{Hom}^{\mrm{cont}}_{C}(M, N) \cong
\opn{lim}_{i \to} \opn{Hom}_{C}(M / \mrm{F}^i M, N) . \]
In this case we consider each
\[ \mrm{F}_i \opn{Hom}^{\mrm{cont}}_{C}(M, N) :=
\opn{Hom}_{C}(M / \mrm{F}^i M, N) \] 
as a discrete inv module, and this endows 
$\opn{Hom}^{\mrm{cont}}_{C}(M, N)$ with a dir-inv structure.

\begin{exa}
In the situation of Example \ref{exa2.4} one has
\[ \opn{Hom}^{\mrm{cont}}_{C}(N, M) \cong
L \otimes_{C} M \]
as dir-inv $C$-modules.
\end{exa}

\begin{exa} \label{exa1.2}
This example is taken from \cite{Ye1}. Assume $\K$ is
noetherian and $C$ a finitely generated $\K$-algebra. For 
$q \in \mbb{N}$ define 
$\mcal{B}_q(C) = \mcal{B}^{-q}(C) := C^{\otimes (q+2)} = 
C \otimes_{\K} \cdots \otimes_{\K} C$.
Define $\what{\mcal{B}}_{q}(C) = \what{\mcal{B}}^{-q}(C)$
to be the adic completion of $\mcal{B}_q(C)$ with respect to the 
ideal $\opn{Ker}(\mcal{B}_q(C) \to C)$.

There is a $\K$-algebra homomorphism 
$\what{\mcal{B}}^{0}(C) \to \what{\mcal{B}}^{-q}(C)$,
corresponding to the two extreme tensor factors, 
and in this way we view $\what{\mcal{B}}^{-q}(C)$ as a
complete inv $\what{\mcal{B}}^{0}(C)$-module. 
There is a continuous coboundary operator that makes
$\what{\mcal{B}}(C) := \boplus_{q \in \mbb{N}}
\what{\mcal{B}}^{-q}(C)$
into a complex of $\what{\mcal{B}}^{0}(C)$-modules, and there is a 
quasi-isomorphism 
$\what{\mcal{B}}(C) \to C$.
We call $\what{\mcal{B}}(C)$ the {\em complete un-normalized bar 
complex} of $C$.

Next define
$\what{\mcal{C}}_{q}(C) = \what{\mcal{C}}^{-q}(C) := 
C \otimes_{\what{\mcal{B}}^{0}(C)} \what{\mcal{B}}^{-q}(C)$.
This is a complete inv $C$-module.
The complex $\what{\mcal{C}}(C)$ is called the 
{\em complete Hochschild chain complex} of $C$.
Finally let
$\mcal{C}_{\mrm{cd}}^{q}(C) := 
\opn{Hom}_{C}^{\mrm{cont}}(\what{\mcal{C}}^{-q}(C), C)$.
The complex
$\mcal{C}_{\mrm{cd}}(C) := \boplus_{q \in \mbb{N}}
\mcal{C}_{\mrm{cd}}^q(C)$ is called the 
{\em continuous Hochschild cochain complex} of $C$.
\end{exa}

\section{Poly Differential Operators}
\label{sec2} 

In this section $\K$ is a commutative base ring, and
$C$ is a commutative $\K$-algebra. 
The symbol $\otimes$ means $\otimes_{\K}$. 
We discuss some basic properties of poly 
differential operators, expanding results from \cite{Ye2}.

\begin{dfn} \label{dfn1.5}
Let $M_1, \ldots, M_p, N$ be $C$-modules. 
A $\K$-multilinear function 
$\phi : M_1 \times \cdots \times M_p \to N$ is called a 
{\em poly differential operator} (over $C$ relative to $\K$)
if there exists some 
$d \in \mbb{N}$ such that for any $(m_1, \ldots, m_p) \in \prod M_i$
and any $i \in \{ 1, \ldots, p \}$ the function 
$M_i \to N$, $m \mapsto \phi(m_1, \ldots, m_{i-1}, m, m_{i+1}, 
\ldots, m_p)$
is a differential operator of order $\leq d$, in the sense of 
\cite[Section 16.8]{EGA-IV}. 
In this case we say that $\phi$ has order $\leq d$ 
in each argument.
\end{dfn}

We shall denote the set of poly differential operators
$\sprod M_i \to N$ 
over $C$ relative to $\K$, of order $\leq d$ in all arguments, 
by
\[ \mrm{F}_d \mcal{D}\mathit{iff}_{\mrm{poly}}
(C; M_1, \ldots, M_p; N) . \] 
And we define
\[ \mcal{D}\mathit{iff}_{\mrm{poly}}
(C; M_1, \ldots, M_p; N) := \bigcup_{d \geq 0}
\mrm{F}_d \mcal{D}\mathit{iff}_{\mrm{poly}}
(C; M_1, \ldots, M_p; N) , \]
the union being inside the set of all $\K$-multilinear functions
$\prod M_i \to N$. 
By default we only consider poly differential operators relative 
to $\K$. 

For a natural number $p$ the $p$-th un-normalized 
$\mcal{B}_p(C)$ was defined in Example \ref{exa1.2}.
Let $I_p(C)$ be the kernel of the ring homomorphism
$\mcal{B}_p(C) \to C$. Define
\[ \mcal{C}_p(C) := C \otimes_{\mcal{B}_0(C)} 
\mcal{B}_p(C) , \]
the {\em $p$-th Hochschild chain module} of $C$ 
(relative to $\K$). For any $d \in \mbb{N}$ define 
\[ \mcal{B}_{p, d}(C) := 
\mcal{B}_{p}(C) \, / \, I_p(C)^{d + 1} , \]
\[ \mcal{C}_{p, d}(C) := C \otimes_{\mcal{B}_{0}(C)}
\mcal{B}_{p, d}(C) \]
and
\[ \mcal{C}_{p, d}(C; M_1, \ldots, M_p) :=
\mcal{C}_{p, d}(C) \otimes_{\mcal{B}_{p-2}(C)}
(M_1 \otimes \cdots \otimes M_p) . \]
Let 
\[ \phi_{\mrm{uni}} : \prod_{i=1}^p M_i \to 
\mcal{C}_{p, d}(C; M_1, \ldots, M_p) \]
be the $\K$-multilinear function 
\[ \phi_{\mrm{uni}}(m_1, \ldots, m_p) := 1 \otimes (m_1 \otimes
\cdots \otimes m_p) . \] 

Observe that for $p = 1$ we get 
$\mcal{C}_{1, d}(C) = \mcal{P}^d (C)$,
the module of principal parts of order $d$ (see \cite{EGA-IV}).
In the same way that 
$\mcal{P}^d (C)$ parameterizes differential operators, 
$\mcal{C}_{p, d}(C)$ 
parameterizes poly differential operators:

\begin{lem} \label{lem3.6}
The assignment $\psi \mapsto \psi \circ \phi_{\mrm{uni}}$
is a bijection
\[ \opn{Hom}_{C}
\big( \mcal{C}_{p, d}(C; M_1, \ldots, M_p), N \big) 
\iso  \mrm{F}_d \mcal{D}\mathit{iff}_{\mrm{poly}}
(C; M_1, \ldots, M_p; N) . \]
\end{lem}

\begin{proof}
The same arguments used in \cite[Section 16.8]{EGA-IV}
also apply here. Cf.\ \cite[Section 1.4]{Ye1}.
\end{proof}

In case $M_1 = \cdots = M_p = N = C$ we see that 
\begin{equation} \label{eqn2.1}
\begin{aligned}
& \mcal{D}\mathit{iff}_{\mrm{poly}}
(C; \underset{p}{\underbrace{C, \ldots, C}}; C) \cong
\lim_{d \to} \opn{Hom}_{C}
\big( \mcal{C}_{p, d}(C), C \big) \\
& \qquad \cong
\opn{Hom}_{C}^{\mrm{cont}}
\big( \what{\mcal{C}}_{p}(C), C \big) =
\mcal{C}^p_{\mrm{cd}}(C) ,
\end{aligned}
\end{equation}
with notation of Example \ref{exa1.2}.

\begin{prop} \label{prop1.9}
Suppose $C$ is a finitely generated $\K$-algebra,
with ideal $\mfrak{c} \subset C$. Let  
$M_1, \ldots, M_p, N$ be $C$-modules, and let
$\phi: \prod M_i \to N$ be a multi differential 
operator over $C$ relative to $\K$. Then $\phi$ is 
continuous for the $\mfrak{c}$-adic dir-inv structures on 
$M_1, \ldots, M_p, N$. 
\end{prop}

\begin{proof}
Suppose $\phi$ has order 
$\leq d$ in each of its arguments, and let 
\[ \psi : \mcal{C}_{p, d}(C; M_1, \ldots, M_p) \to N \]
be the corresponding $C$-linear homomorphism.
As in \cite[Proposition 1.4.3]{Ye1}, since $C$ is a finitely 
generated $\K$-algebra, it follows that
$\mcal{B}_{p, d}(C)$ is a finitely generated module over
$\mcal{B}_{0}(C)$; and hence $\mcal{C}_{p, d}(C)$
is a finitely generated $C$-module.
Let's denote by $\set{\mrm{F}_j M_i}_{j \in J_i}$ and 
$\set{\mrm{F}_k N}_{k \in K}$
the $\mfrak{c}$-adic dir-inv structures on $M_i$ and 
$N$. For any $j_1, \ldots, j_p$ the 
$\mcal{B}_{p-2}(C)$-module
$\mrm{F}_{j_1} M_1 \otimes \cdots \otimes \mrm{F}_{j_p} M_p$
is finitely generated, and hence the $C$-module
$\mcal{C}_{p, d}(C; \mrm{F}_{j_1} M_1, \ldots, \mrm{F}_{j_p} 
M_p)$ 
is finitely generated. Therefore
\[ \psi \big( \mcal{C}_{p, d}(C; 
\mrm{F}_{j_1} M_1, \ldots, \mrm{F}_{j_p} M_p) \big) = 
\mrm{F}_k N \]
for some $k \in K$. 

It remains to prove that 
$\phi : \prod_{i=1}^p \mrm{F}_{j_i} M_i \to \mrm{F}_k N$ 
is continuous for the $\mfrak{c}$-adic inv structures. 
But just like \cite[Proposition 1.4.6]{Ye1}, for any $i$ 
and $l$ one has
\begin{equation} \label{eqn3.4}
\phi(\mrm{F}_{j_1} M_1, \ldots,
\mfrak{c}^{i+d} \mrm{F}_{j_l} M_l , \ldots,
\mrm{F}_{j_p} M_p) \subset 
\mfrak{c}^{i} \mrm{F}_k N .
\end{equation}
\end{proof} 

Suppose $C'$ is a commutative $C$-algebra with ideal 
$\mfrak{c}' \subset C'$. One says that $C'$ is 
{\em $\mfrak{c}'$-adically 
formally \'etale} over $C$ if the following condition holds. Let 
$D$ be a commutative $C$-algebra with nilpotent ideal
$\mfrak{d}$, and let $f : C' \to D / \mfrak{d}$ be a $C$-algebra 
homomorphism such that $f({\mfrak{c}'}^i) = 0$ for 
$i \gg 0$. Then $f$ lifts uniquely to a $C$-algebra homomorphism 
$\til{f} : C' \to D$. The important instances are when $C \to C'$ 
is \'etale (and then $\mfrak{c}' = 0$); and when 
$C'$ is the $\mfrak{c}$-adic completion of $C$ for some ideal 
$\mfrak{c} \subset A$ (and $\mfrak{c}' = C' \mfrak{c}$).
In both these instances $C'$ is $\mfrak{c}$-adically complete; 
and if $C$ is noetherian, then $C \to C'$ is also flat.

\begin{lem} \label{lem3.7}
Let $C'$ be a $\mfrak{c}'$-adically formally \'etale $C$-algebra. 
Define $C'_j := C' / {\mfrak{c}'}^{j+1}$. Consider $C'$ and 
$\mcal{C}_{p, d}(C)$ as inv $C$-modules, with the 
$\mfrak{c}'$-adic and discrete inv structures respectively.
Then the canonical homomorphism
\[ C' \hatotimes{C} \mcal{C}_{p, d}(C) \to
\opn{lim}_{\leftarrow j} \mcal{C}_{p, d}(C'_j) \]
is bijective.
\end{lem}

\begin{proof}
Define ideals
\[ \mfrak{c}'_p :=  \opn{Ker} \big( \mcal{C}_{p}(C')
\to \mcal{C}_{p}(C'_0) \big)  \] 
and
\[ J := \opn{Ker} \big( C'_j \otimes_C \mcal{C}_{p, d}(C)
\to C'_j \big) . \] 
By the transitivity and the base change properties of formally
\'etale homomorphisms, the ring homomorphism
\[ \mcal{C}_{p}(C) \cong C \otimes \cdots \otimes C
\to C' \otimes \cdots \otimes C' \cong \mcal{C}_{p}(C') 
\]
is $\mfrak{c}_p'$-adically formally \'etale.
Consider the commutative diagram of ring homomorphisms
(with solid arrows)
\[ \UseTips \xymatrix @=8ex {
C \ar[r] \ar[d]
& \mcal{C}_{p}(C) \ar[r] \ar[d]
& C'_j \otimes_C \mcal{C}_{p, d}(C) \ar[r]^{e} \ar[d]
& \mcal{C}_{p, d}(C'_j) \ar[d] \\
C' \ar[r] \ar@{-->}[urr]
& \mcal{C}_{p}(C') \ar[r]_{f} \ar@{-->}[ur]|{\til{f}}
\ar@{-->}[urr]_(0.7){g}
&  C'_j \ar[r]^{=}
& C'_j \ .
} \]
The ideal $J$ satisfies $J^{d+1} = 0$, and the ideal
$\opn{Ker} \big( \mcal{C}_{p, d}(C'_j) \to C'_j \big)$
is nilpotent too. 
Due to the unique lifting property the dashed arrows exist and are 
unique, making the whole diagram commutative. Moreover 
$g : \mcal{C}_{p}(C') \to \mcal{C}_{p}(C'_j)$
has to be the canonical surjection, and $\til{f}$ is surjective. 

A little calculation shows that 
$\til{f}(I_p(C')^{d+1}) = 0$, 
and hence $\til{f}$ induces a homomorphism
\[ \bar{f} : \mcal{C}_{p, d}(C') \to
C'_j \otimes_C \mcal{C}_{p, d}(C) . \] 
Let
\[ \mfrak{c}_{p, d}' := \opn{Ker} \big( \mcal{C}_{p, d}(C') \to
\mcal{C}_{p, d}(C'_0) \big) . \]
Another calculation shows that 
$\bar{f}(\mfrak{c}_{p, d}' {}^{(j+1)(d+1)}) = 0$.
The conclusion is that there are surjections
\[ \mcal{C}_{p, d}(C'_{jd + j + d}) \xar{\bar{f}}
C'_j \otimes_C \mcal{C}_{p, d}(C) \xar{e}
\mcal{C}_{p, d}(C'_j) , \]
such that $e \circ \bar{f}$ is the canonical surjection.
Passing to the inverse limit we deduce that 
\[ C' \hatotimes{C} \mcal{C}_{p, d}(C) \to
\opn{lim}_{\leftarrow j} \mcal{C}_{p, d}(C'_j) \]
is bijective.
\end{proof}

\begin{prop} \label{prop3.9}
Assume $C$ is a noetherian finitely generated $\K$-algebra, and 
$C'$ is a noetherian, $\mfrak{c}'$-adically complete, flat,
$\mfrak{c}'$-adically formally \'etale $C$-algebra. 
Let $M_1, \ldots, M_p, N$ be $C$-modules, and define
$M'_i := C' \otimes_C M_i$ and $N' := C' \otimes_C N$.
\begin{enumerate}
\item Suppose $\phi: \prod_{i=1}^p M_i \to N$ 
is a poly differential operator over $C$. 
Then $\phi$ extends uniquely 
to a poly differential operator 
$\phi' : \prod_{i=1}^p M_i' \to N'$ over $C'$. 
If $\phi$ has order $\leq d$ then so does $\phi'$.
\item The homomorphism
\[ \begin{aligned}
& C' \otimes_C \mrm{F}_d \mcal{D}\mathit{iff}_{\mrm{poly}}
(C; M_1, \ldots, M_p; N) \\
& \qquad \to
\mrm{F}_d \mcal{D}\mathit{iff}_{\mrm{poly}}
(C'; M'_1, \ldots, M'_p; N') , 
\end{aligned} \]
$c' \otimes \phi \mapsto c' \phi$, 
is bijective.
\end{enumerate}
\end{prop}

\begin{proof}
By Proposition \ref{prop1.9}, applied to $C$ with the $0$-adic inv 
structure, we may assume that the $C$-modules 
$M_1, \ldots, M_p, N$ are finitely generated. 

Fix $d \in \mbb{N}$. Define 
$C'_j := C' / \mfrak{c}'{}^{j+1}$ and 
$N'_j := C'_j \otimes_C N$. 
So 
$C' \cong \lim_{\leftarrow j} C'_j$
and $N' \cong \lim_{\leftarrow j} N'_j$.

By Lemma \ref{lem3.6} and Proposition \ref{prop1.9} we have
\begin{equation} \label{eqn3.3}
\begin{aligned}
& \mrm{F}_d \mcal{D}\mathit{iff}_{\mrm{poly}}(C'; 
M'_1, \ldots, M'_p; N') \\
& \qquad \cong
\opn{Hom}_{C'} \big( \mcal{C}_{p, d}(C'; 
M'_1, \ldots, M'_p), N' \big) \\
& \qquad \cong 
\lim_{\leftarrow j} 
\opn{Hom}_{C'} \big( \mcal{C}_{p, d}(C';
M'_1, \ldots, M'_p), N'_j \big) . 
\end{aligned}
\end{equation}
Now for any $k \geq j + d$ one has
\[ \opn{Hom}_{C'} \big( \mcal{C}_{p, d}(C';
M'_1, \ldots, M'_p), N'_j \big)
\cong
\opn{Hom}_{C'} \big( \mcal{C}_{p, d}(C'_k;
M'_1, \ldots, M'_p), N'_j \big) . \]
This is because of formula (\ref{eqn3.4}). 
Thus, using Lemma \ref{lem3.7}, we obtain
\[ \begin{aligned}
& \opn{Hom}_{C'} \big( \mcal{C}_{p, d}(C';
M'_1, \ldots, M'_p), N'_j \big) \\
& \qquad \cong
\opn{Hom}_{C'} \big( \lim_{\leftarrow k} 
\mcal{C}_{p, d}(C'_k; M'_1, \ldots, M'_p), N'_j \big) \\
& \qquad \cong 
\opn{Hom}_{C'} \big( C' \otimes_C 
\mcal{C}_{p, d}(C; M_1, \ldots, M_p), N'_j \big) \\
& \qquad \cong 
\opn{Hom}_{C} \big( \mcal{C}_{p, d}(C; M_1, \ldots, M_p)
, N'_j \big) . 
\end{aligned} \]
Combining this with (\ref{eqn3.3}) we get
\[ \begin{aligned}
& \mrm{F}_d \mcal{D}\mathit{iff}_{\mrm{poly}}(C'; 
M'_1, \ldots, M'_p; N') \\
& \qquad \cong
\lim_{\leftarrow j} 
\opn{Hom}_{C} \big( \mcal{C}_{p, d}(C; M_1, \ldots, M_p), 
N'_j \big) \\
& \qquad \cong 
\opn{Hom}_{C} \big( \mcal{C}_{p, d}(C; M_1, \ldots, M_p), 
N' \big)  . 
\end{aligned} \]
But $C \to C'$ is flat, $C$ is noetherian, and 
$\mcal{C}_{p, d}(C; M_1, \ldots, M_p)$ 
is a finitely generated $C$-module. 
Therefore 
\[ \begin{aligned}
& \opn{Hom}_{C} \big( \mcal{C}_{p, d}(C; M_1, \ldots, M_p), 
N' \big) \\
& \qquad \cong
C' \otimes_C 
\opn{Hom}_{C} \big( \mcal{C}_{p, d}(C; M_1, \ldots, M_p), 
N \big)  . 
\end{aligned} \]
The conclusion is that
\begin{equation} \label{eqn3.5}
\begin{aligned}
& \mrm{F}_m \mcal{D}^{p+1}_{\mrm{poly}}(C'; 
M'_1, \ldots, M'_p; N') \\
& \qquad \cong
C' \otimes_C \mrm{F}_m \mcal{D}^{p+1}_{\mrm{poly}}(C;
M_1, \ldots, M_p; N) .
\end{aligned} 
\end{equation} 

Given $\phi : \sprod M_i \to N$ of order $\leq d$, let
$\phi' := 1 \otimes \phi$ under the isomorphism (\ref{eqn3.5}).
Backtracking we see that $\phi'$ is the unique poly differential 
operator extending $\phi$. 
\end{proof}

\section{$\mrm{L}_{\infty}$ Morphisms and their Twists}
\label{sec3}

In this section we expand some results on $\mrm{L}_{\infty}$ 
algebras and morphisms from \cite{Ko1} Section 4. Much of the 
material presented here is based on discussions with
Vladimir Hinich. There is some overlap with Section 2.2 
of \cite{Fu}, with Section 6.1 of \cite{Le}, and possibly with 
other accounts. 

Let $\K$ be a field of characteristic $0$.
Given a graded $\K$-module
$\mfrak{g} = \boplus_{i \in \mbb{Z}} \mfrak{g}^i$ 
and a natural number $j$ let 
$\mrm{T}^j \mfrak{g} :=$
$\underset{j}{\underbrace{\mfrak{g} \otimes \cdots \otimes
\mfrak{g}}}$. 
The direct sum
$\mrm{T} \mfrak{g} := \boplus_{j \in \mbb{N}} \mrm{T}^j \mfrak{g}$
is the tensor algebra. 
Let us denote the multiplication in $\mrm{T} \mfrak{g}$ by 
$\circledast$. (This is just another way of writing $\otimes$, but 
it will be convenient to do so.) 

The permutation group $\mfrak{S}_j$ acts on 
$\mrm{T}^j \mfrak{g}$ as follows. For any sequence of integers
$\bsym{d} = (d_1, \ldots, d_j)$ there is a group homomorphism
$\opn{sgn}_{\bsym{d}} : \mfrak{S}_j \to \{ \pm1 \}$
such that on a transposition $\sigma = (p, p+1)$ the value is
$\opn{sgn}_{\bsym{d}}(\sigma) = (-1)^{d_p d_{p+1}}$.
The action of a permutation
$\sigma \in \mfrak{S}_j$ on $\mrm{T}^j \mfrak{g}$ 
is then
\[ \sigma (\gamma_1 \circledast \cdots \circledast \gamma_j) :=
\opn{sgn}_{\bsym{d}}(\sigma)
\gamma_{\sigma(1)} \circledast \cdots \circledast 
\gamma_{\sigma(j)} \]
for 
$\gamma_1 \in \mfrak{g}^{d_1}, \ldots, 
\gamma_j \in \mfrak{g}^{d_j}$.
Define $\til{\mrm{S}}^j \mfrak{g}$ to be the set of 
$\mfrak{S}_j$-invariants inside $\mrm{T}^j \mfrak{g}$, and
$\til{\mrm{S}} \mfrak{g} := 
\boplus_{j \geq 0} \til{\mrm{S}}^j \mfrak{g}$.

The $\K$-module $\mrm{T} \mfrak{g}$ is also a coalgebra, with 
coproduct
$\til{\Delta} : \mrm{T} \mfrak{g} \to
\mrm{T} \mfrak{g} \otimes \mrm{T} \mfrak{g}$
given by the formula
\[ \til{\Delta}(\gamma_1 \circledast \cdots \circledast \gamma_j)
:= \sum_{p = 0}^j 
(\gamma_1 \circledast \cdots \circledast \gamma_p) \otimes
(\gamma_{p+1} \circledast \cdots \circledast \gamma_j) . \]
The submodule 
$\til{\mrm{S}} \mfrak{g} \subset \mrm{T} \mfrak{g}$
is a sub-coalgebra (but not a subalgebra!). 

The super-symmetric algebra
$\mrm{S} \mfrak{g} = \boplus_{j \geq 0} \mrm{S}^j \mfrak{g}$
is defined to be the 
quotient of $\mrm{T} \mfrak{g}$ by the ideal generated by the 
elements 
$\gamma_1 \circledast \gamma_2 - (-1)^{d_1 d_2}
\gamma_2 \circledast \gamma_1$,
for all $\gamma_1 \in \mfrak{g}^{d_1}$ and
$\gamma_2 \in \mfrak{g}^{d_2}$. 
In other words, $\mrm{S}^j \mfrak{g}$ is the set of coinvariants 
of $\mrm{T}^j \mfrak{g}$ under the action of the group 
$\mfrak{S}_j$. The product in the algebra $\mrm{S} \mfrak{g}$ is 
denoted by $\cdot$. The canonical projection is
$\pi : \mrm{T} \mfrak{g} \to \mrm{S} \mfrak{g}$
is an algebra homomorphism:
$\pi(\gamma_1 \circledast \gamma_2) = \gamma_1 \cdot
\gamma_2$. 

In fact $\mrm{S} \mfrak{g}$ is a commutative 
cocommutative Hopf algebra. The comultiplication 
\[ \Delta : \mrm{S} \mfrak{g} \to \mrm{S} \mfrak{g} \otimes
\mrm{S} \mfrak{g} \]
is the unique $\K$-algebra homomorphism such that
\[ \Delta(\gamma) = \gamma \otimes 1 + 1 \otimes \gamma \]
for all $\gamma \in \mfrak{g}$. 
The antipode is $\gamma \mapsto - \gamma$. 
The projection $\pi : \mrm{T} \mfrak{g} \to \mrm{S} \mfrak{g}$ 
is not a coalgebra homomorphism. However:

\begin{lem} \label{lem3.3}
Let $\tau : \mrm{S} \mfrak{g} \to \mrm{T} \mfrak{g}$ be the 
$\K$-module homomorphism defined by
\[ \tau(\gamma_1 \cdots  \gamma_j) :=
\sum_{\sigma \in \mfrak{S}_j} 
\opn{sgn}_{(d_1, \ldots, d_j)}(\sigma)
\gamma_{\sigma(1)} \circledast \cdots \circledast 
\gamma_{\sigma(j)} \]
for 
$\gamma_1 \in \mfrak{g}^{d_1}, \ldots, 
\gamma_j \in \mfrak{g}^{d_j}$.
Then $\tau : \mrm{S} \mfrak{g} \to \til{\mrm{S}} \mfrak{g}$
is a coalgebra isomorphism, where $\mrm{S} \mfrak{g}$ has the 
comultiplication $\Delta$ and $\til{\mrm{S}} \mfrak{g}$ has the 
comultiplication $\til{\Delta}$.
\end{lem}

\begin{proof}
Define $\til{\pi} : \mrm{T} \mfrak{g} \to \mrm{S} \mfrak{g}$
to be the $\K$-module homomorphism
\[ \til{\pi}(\gamma_1 \circledast \cdots \circledast \gamma_j)
:= \smfrac{1}{j!} 
\pi(\gamma_1 \circledast \cdots \circledast \gamma_j) = 
\smfrac{1}{j!} \gamma_1 \cdots  \gamma_j  \]
for $\gamma_1, \ldots, \gamma_j \in \mfrak{g}$. 
So 
$\til{\pi} \circ \tau$ is the identity map of $\mrm{S} \mfrak{g}$,
and $\til{\pi} : \til{\mrm{S}} \mfrak{g} \to
\mrm{S} \mfrak{g}$
is bijective. It suffices to prove that
\[ (\til{\pi} \otimes \til{\pi}) \circ (\tau \otimes \tau) \circ 
\Delta = 
(\til{\pi} \otimes \til{\pi}) \circ \til{\Delta} \circ \tau . \]

Take any 
$\gamma_1 \in \mfrak{g}^{d_1}, \ldots, 
\gamma_j \in \mfrak{g}^{d_j}$ and write
$\bsym{d} := (d_1, \ldots, d_j)$. 
Then
\[ \begin{aligned}
& \bigl( (\til{\pi} \otimes \til{\pi}) \circ \til{\Delta} \circ 
\tau \bigr) (\gamma_1 \cdots  \gamma_j) \\
& \qquad = 
\sum_{p =0}^j \ \sum_{\sigma \in \mfrak{S}_j}
\smfrac{1}{p!} \smfrac{1}{(j-p)!} \opn{sgn}_{\bsym{d}}(\sigma)
(\gamma_{\sigma(1)} \cdots \gamma_{\sigma(p)}) \otimes
(\gamma_{\sigma(p+1)} \cdots \gamma_{\sigma(j)}) .
\end{aligned} \]
On the other hand
\[ \begin{aligned}
& \bigl( (\til{\pi} \otimes \til{\pi}) \circ (\tau \otimes \tau)
\circ \Delta \bigr) (\gamma_1 \cdots  \gamma_j) \\
& \qquad = 
\Delta(\gamma_1 \cdots  \gamma_j) = 
(1 \otimes \gamma_1 + \gamma_1 \otimes 1) \cdots 
(1 \otimes \gamma_j + \gamma_j \otimes 1) \\
& \qquad 
\sum_{p =0}^j \ \sum_{\sigma \in \mfrak{S}_{p, j-p}}
\opn{sgn}_{\bsym{d}}(\sigma)
(\gamma_{\sigma(1)} \cdots \gamma_{\sigma(p)}) \otimes
(\gamma_{\sigma(p+1)} \cdots \gamma_{\sigma(j)}) ,
\end{aligned} \]
where $\mfrak{S}_{p, j-p}$ is the set of $(p, j-p)$-shuffles 
inside the group $\mfrak{S}_{j}$. 
Since the algebra $\mrm{S} \mfrak{g}$ is super-commutative the two 
sums are equal.
\end{proof}

The grading on $\mfrak{g}$ induces a grading on 
$\mrm{S} \mfrak{g}$, which we call the {\em degree}. Thus for 
$\gamma_i \in \mfrak{g}^{d_i}$ the degree of 
$\gamma_1 \cdots \gamma_j \in \mrm{S}^j \mfrak{g}$ 
is $d_1 + \cdots + d_j$
(unless $\gamma_1 \cdots \gamma_j = 0$). We consider 
$\mrm{S} \mfrak{g}$ as a graded algebra for this grading. 
Actually there is another grading on $\mrm{S} \mfrak{g}$, by 
{\em order}, where we define the order of 
$\gamma_1 \cdots \gamma_j$ to be
$j$ (again, unless this element is zero). But this grading will 
have a different role.
 
By definition the $j$-th super-exterior power of $\mfrak{g}$ is
\begin{equation} \label{eqn10.6}
\bwedge^{j} \mfrak{g} :=
\mrm{S}^j(\mfrak{g}[1])[-j] ,
\end{equation}
where $\mfrak{g}[1]$ is the shifted graded module whose degree $i$ 
component is $\mfrak{g}[1]^i = \mfrak{g}^{i+1}$. 
When $\mfrak{g}$ is concentrated 
in degree $0$ then these are the usual 
constructions of symmetric and exterior algebras, respectively. 

We denote by 
$\opn{ln} : \mrm{S} \mfrak{g} \to \mrm{S}^1 \mfrak{g}
= \mfrak{g}$
the projection. So $\opn{ln}(\gamma)$ is the $1$-st order term of 
$\gamma \in \mrm{S} \mfrak{g}$. 
(The expression ``$\opn{ln}$'' might stand for ``linear'' or 
``logarithm''.)

\begin{dfn}
Let $\mfrak{g}$ and $\mfrak{g}'$ be two graded $\K$-modules, and 
let $\Psi : \mrm{S} \mfrak{g} \to \mrm{S} \mfrak{g}'$
be a $\K$-linear homomorphism. For any $j \geq 1$ the
{\em $j$-th Taylor coefficient} of $\Psi$ is defined to be
\[ \partial^j \Psi := \opn{ln} \circ \Psi :
\mrm{S}^j \mfrak{g} \to \mfrak{g}' . \]
\end{dfn}

\begin{lem}
Suppose we are given a sequence of $\K$-linear homomorphisms 
$\psi_j : \mrm{S}^j \mfrak{g} \to \mfrak{g}'$,
$j \geq 1$, each of degree $0$. Then there is a unique 
coalgebra homomorphism
$\Psi : \mrm{S} \mfrak{g} \to \mrm{S} \mfrak{g}'$,
homogeneous of degree $0$ and satisfying $\Psi(1) = 1$,
whose Taylor coefficients are
$\partial^j \Psi =\psi_j$.
\end{lem}

\begin{proof}
Let 
$\til{\opn{ln}} : \til{\mrm{S}} \mfrak{g}' \to 
\til{\mrm{S}}^1 \mfrak{g}' = \mfrak{g}'$
be the projection for this coalgebra. 
Consider the exact sequence of coalgebras
\begin{equation} \label{eqn3.6}
0 \to \K \to \til{\mrm{S}} \mfrak{g} \to
\til{\mrm{S}}^{\geq 1} \mfrak{g} \to 0 . 
\end{equation}
According to \cite[Section 4.1]{Ko1} (see also 
\cite[Lemma 2.1.5]{Fu}) the sequence 
$\{ \psi_j \}_{j \geq 1}$
uniquely determines a coalgebra homomorphism
$\til{\Psi} : \til{\mrm{S}}^{\geq 1} \mfrak{g} \to 
\til{\mrm{S}}^{\geq 1} \mfrak{g}'$
such that 
\[ \til{\opn{ln}} \circ \til{\Psi} |_{\til{\mrm{S}}^j \mfrak{g}}
= \psi_j \circ \tau^{-1} |_{\til{\mrm{S}}^j \mfrak{g}} \]
for all $j \geq 1$. Here
$\tau : \mrm{S} \mfrak{g} \iso \til{\mrm{S}} \mfrak{g}$
is the coalgebra isomorphism of Lemma \ref{lem3.3}.
Using (\ref{eqn3.6}) we can lift $\til{\Psi}$ uniquely to a 
coalgebra homomorphism 
$\til{\Psi} : \til{\mrm{S}} \mfrak{g} \to 
\til{\mrm{S}} \mfrak{g}'$
by setting $\til{\Psi}(1) := 1$. Now
define the coalgebra homomorphism 
$\Psi : \mrm{S} \mfrak{g} \to \mrm{S} \mfrak{g}'$
to be
$\Psi := \tau^{-1} \circ \til{\Psi} \circ \tau$.
\end{proof}

A $\K$-linear map $Q : \mrm{S} \mfrak{g} \to \mrm{S} \mfrak{g}$ 
is a {\em coderivation} if 
\[ \Delta \circ Q = (Q \otimes \bsym{1} + \bsym{1} \otimes Q)
\circ \Delta , \]
where $\bsym{1} := \bsym{1}_{\mrm{S} \mfrak{g}}$, the identity 
map. 

\begin{lem}
Given a sequence of $\K$-linear homomorphisms 
$\psi_j : \mrm{S}^j \mfrak{g} \to \mfrak{g}$,
$j \geq 1$, each of degree $1$, there is a unique 
coderivation $Q$ of degree $1$, such that  $Q(1) = 0$ and 
$\partial^j Q = \psi_j$. Furthermore, one has
$Q(\gamma) \in \boplus_{j \geq 1} \mrm{S}^j \mfrak{g}$
for every $\gamma \in \mrm{S} \mfrak{g}$.
\end{lem}

\begin{proof}
According to \cite[Section 4.3]{Ko1} (see also 
\cite[Lemma 2.1.2]{Fu}) the sequence 
$\{ \psi_j \}_{j \geq 1}$
uniquely determines a coderivation 
$\til{Q} : \til{\mrm{S}}^{\geq 1} \mfrak{g} \to 
\til{\mrm{S}}^{\geq 1} \mfrak{g}$
such that 
\[ \til{\opn{ln}} \circ \til{Q} |_{\til{\mrm{S}}^j \mfrak{g}}
= \psi_j \circ \tau^{-1} |_{\til{\mrm{S}}^j \mfrak{g}} \]
for all $j \geq 1$. Using (\ref{eqn3.6}) this can be lifted
uniquely to a coderivation 
$\til{Q} : \til{\mrm{S}} \mfrak{g} \to 
\til{\mrm{S}} \mfrak{g}$
by setting $\til{Q}(1) := 0$. Now
define the coderivation
$Q : \mrm{S} \mfrak{g} \to \mrm{S} \mfrak{g}$
to be
$Q := \tau^{-1} \circ \til{Q} \circ \tau$. 
\end{proof}

We will be mostly interested in the coalgebras 
$\mrm{S}(\mfrak{g}[1])$ and $\mrm{S}(\mfrak{g}'[1])$.
Observe that if 
$\Psi : \mrm{S}(\mfrak{g}[1]) \to \mrm{S}(\mfrak{g}'[1])$
is a homogeneous $\K$-linear homomorphism of degree $i$, then,
using formula (\ref{eqn10.6}), each Taylor coefficient 
$\partial^j \Psi$ may be viewed as a homogeneous $\K$-linear 
homomorphism
$\partial^j \Psi : \bwedge^j \mfrak{g} \to \mfrak{g}$
of degree $i + 1 - j$. 

\begin{dfn}
Let $\mfrak{g}$ be a graded $\K$-module. An  $\mrm{L}_{\infty}$ 
algebra structure on $\mfrak{g}$ is a coderivation
$Q : \mrm{S}(\mfrak{g}[1]) \to \mrm{S}(\mfrak{g}[1])$ 
of degree $1$, satisfying $Q(1) = 0$ and 
$Q \circ Q = 0$. We call the pair 
$(\mfrak{g}, Q)$ an {\em $\mrm{L}_{\infty}$ algebra}.
\end{dfn}

The notion of $\mrm{L}_{\infty}$ algebra 
generalizes that of DG Lie algebra in the following sense:

\begin{prop}[{\cite[Section 4.3]{Ko1}}] \label{prop9.4}
Let $Q : \mrm{S}(\mfrak{g}[1]) \to \mrm{S}(\mfrak{g}[1])$ be a 
coderivation of degree $1$ with $Q(1) = 0$. 
Then the following conditions are equivalent.
\begin{enumerate}
\rmitem{i} $\partial^j Q = 0$ for all $j \geq 3$, and 
$Q \circ Q = 0$.
\rmitem{ii} $\partial^j Q = 0$ for all $j \geq 3$, and
$\mfrak{g}$ is a DG Lie algebra with respect to the
differential $\mrm{d} := \partial^1 Q$ and the bracket 
$[-,-] := \partial^2 Q$.
\end{enumerate}
\end{prop}

In view of this, we shall say that $(\mfrak{g}, Q)$ is a DG Lie 
algebra if the equivalent conditions of the proposition hold.
An easy calculation shows that given an 
$\mrm{L}_{\infty}$ algebra $(\mfrak{g}, Q)$,
the function $\partial^1 Q : \mfrak{g} \to \mfrak{g}$ 
is a differential, and
$\partial^2 Q$ induces a graded Lie bracket on 
$\mrm{H} (\mfrak{g}, \partial^1 Q)$. We shall denote this graded 
Lie algebra by $\mrm{H}(\mfrak{g}, Q)$.

\begin{dfn}
Let $(\mfrak{g}, Q)$ and $(\mfrak{g}', Q')$ be $\mrm{L}_{\infty}$
algebras. An {\em $\mrm{L}_{\infty}$ morphism  
$\Psi : (\mfrak{g}, Q) \to (\mfrak{g}', Q')$} 
is a coalgebra homomorphism
$\Psi : \mrm{S}(\mfrak{g}[1]) \to \mrm{S}(\mfrak{g}'[1])$ 
of degree $0$, satisfying $\Psi(1) = 1$ and 
$\Psi \circ Q = Q' \circ \Psi$.
\end{dfn}

\begin{prop}[{\cite[Section 4.3]{Ko1}}] \label{prop9.1}
Let $(\mfrak{g}, Q)$ and $(\mfrak{g}', Q')$ be DG Lie algebras,
and let $\Psi : \mrm{S}(\mfrak{g}[1]) \to \mrm{S}(\mfrak{g}'[1])$ 
be a coalgebra homomorphism of degree $0$ such that $\Psi(1) = 1$. 
Then $\Psi$ is an $\mrm{L}_{\infty}$ morphism \tup{(}i.e.\ 
$\Psi \circ Q = Q' \circ \Psi$\tup{)} iff the Taylor coefficients
$\psi_i := \partial^i \Psi : \bwedge^i \mfrak{g} \to \mfrak{g}'$
satisfy the following identity:
\[ \begin{aligned}
& \mrm{d} \bigl( \psi_i(\gamma_1 \wedge \cdots 
\wedge \gamma_i) \bigr) - \sum_{k = 1}^i \pm 
\psi_i \bigl( \gamma_1 \wedge \cdots \wedge 
\mrm{d}(\gamma_k) \wedge \cdots \wedge \gamma_i \bigr) = \\
& \quad \smfrac{1}{2} \sum_{\substack{k, l \geq 1 \\ k + l = i}}
\smfrac{1}{k! \, l!} \sum_{\sigma \in \mfrak{S}_i} \pm
\bigl[ \psi_k (\gamma_{\sigma(1)} \wedge \cdots \wedge 
\gamma_{\sigma(k)}),
\psi_l (\gamma_{\sigma(k + 1)} \wedge \cdots \wedge 
\gamma_{\sigma(i)}) \bigr] \\
& \qquad +
\sum_{k < l} \pm
\psi_{i-1} \bigl( [\gamma_k, \gamma_l] \wedge
\gamma_{1} \wedge \cdots \gamma_k \hspace{-1em} \diagup \cdots
\gamma_l \hspace{-1em} \diagup 
\cdots \wedge \gamma_{i} \bigr) .
\end{aligned} \]
Here $\gamma_k \in \mfrak{g}$ are homogeneous elements,
$\mfrak{S}_i$ is the permutation group of $\{ 1, \ldots, i \}$,
and the signs depend only on the indices, the permutations and the 
degrees of the elements $\gamma_k$.
\tup{(}See \cite[Section 6]{Ke} or
\cite[Theorem 3.1]{CFT} for the explicit signs.\tup{)}
\end{prop}

The proposition shows that when $(\mfrak{g}, Q)$ and 
$(\mfrak{g}', Q')$ are DG Lie algebras and $\partial^j \Psi = 0$ 
for all $j \geq 2$, then 
$\partial^1 \Psi : \mfrak{g} \to \mfrak{g}'$ 
is a homomorphism of DG Lie algebras; and conversely. It also 
implies that for any $\mrm{L}_{\infty}$ morphism
$\Psi : (\mfrak{g}, Q) \to (\mfrak{g}', Q')$, 
the map
$\mrm{H}(\Psi) : \mrm{H} (\mfrak{g}, Q) \to 
\mrm{H} (\mfrak{g}', Q')$
is a homomorphism of graded Lie algebras. 

Given DG Lie algebras $\mfrak{g}$ and $\mfrak{g}'$ we consider 
them as $\mrm{L}_{\infty}$ algebras
$(\mfrak{g}, Q)$ and $(\mfrak{g}', Q')$, as explained in 
Proposition \ref{prop9.4}. If 
$\Psi : (\mfrak{g}, Q) \to (\mfrak{g}', Q')$
is an $\mrm{L}_{\infty}$ morphism, then we shall say (by slight 
abuse of notation) that
$\Psi : \mfrak{g} \to \mfrak{g}'$
is an $\mrm{L}_{\infty}$ morphism.

 From here until Theorem \ref{thm9.1} (inclusive)
$C$ is a commutative $\K$-algebra, 
and $\mfrak{g}$, $\mfrak{g}'$ are graded $C$-modules. 
Suppose $(\mfrak{g}, Q)$ is an $\mrm{L}_{\infty}$ 
algebra structure on $\mfrak{g}$ such that the Taylor coefficients
$\partial^j Q : \bwedge^j \mfrak{g} \to \mfrak{g}$
are all $C$-multilinear. Then we say $(\mfrak{g}, Q)$ is a
{\em $C$-multilinear $\mrm{L}_{\infty}$ algebra}. Similarly one 
defines the notion of 
{\em $C$-multilinear $\mrm{L}_{\infty}$ morphism}
$\Psi : (\mfrak{g}, Q) \to (\mfrak{g}', Q')$. 

With $C$ and $\mfrak{g}$ as above let 
$\mrm{S}_C \mfrak{g}$ be the super-symmetric associative unital 
free algebra over $C$. Namely $\mrm{S}_C \mfrak{g}$ is the 
quotient of the tensor algebra 
$\mrm{T}_C \mfrak{g} = C \oplus \mfrak{g} \oplus 
(\mfrak{g} \otimes_C \mfrak{g}) \oplus \cdots$
by the ideal generated by the super-commutativity relations. 
The algebra $\mrm{S}_C \mfrak{g}$ is a
Hopf algebra over $C$, with comultiplication 
\[ \Delta_C : \mrm{S}_C \mfrak{g} \to
\mrm{S}_C \mfrak{g} \otimes_C
\mrm{S}_C \mfrak{g} . \]
The formulas are just as in the case $C = \K$. 
It will be useful to note that $\Delta_C$ preserves the grading by 
order, namely
\[ \Delta_C(\mrm{S}_C^i \mfrak{g}) \subset 
\boplus_{j+k=i} \mrm{S}_C^j \mfrak{g} \otimes_C
\mrm{S}_C^k \mfrak{g} . \]

\begin{lem} \label{lem9.1}
\begin{enumerate}
\item Let $\mfrak{g}$ be a graded $C$-module. There is a canonical 
bijection $Q \mapsto Q_C$
between the set of $C$-multilinear $\mrm{L}_{\infty}$
algebra structures $Q$ on $\mfrak{g}$, and the set of coderivations 
$Q_C : \mrm{S}_C(\mfrak{g}[1]) \to \mrm{S}_C(\mfrak{g}[1])$
over $C$ of degree $1$, such that $Q_C(1) = 0$ and 
$Q_C \circ Q_C = 0$. 
\item Let $(\mfrak{g}, Q)$ and $(\mfrak{g}', Q')$ be two 
$C$-multilinear $\mrm{L}_{\infty}$ algebras. The set of 
$C$-multilinear $\mrm{L}_{\infty}$ morphisms 
$\Psi : (\mfrak{g}, Q) \to (\mfrak{g}', Q')$
is canonically bijective to the set of coalgebra homomorphisms
$\Psi_C : \mrm{S}_C(\mfrak{g}[1]) \to \mrm{S}_C(\mfrak{g}'[1])$
over $C$ of degree $0$, such that $\Psi_C(1) = 1$ and 
$\Psi_C \circ Q_C = Q'_C \circ \Psi_C$. 
\end{enumerate}
\end{lem}

\begin{proof}
The data for a coderivation
$Q_C : \mrm{S}_C(\mfrak{g}[1]) \to \mrm{S}_C(\mfrak{g}[1])$
over $C$ is its sequence of $C$-linear Taylor coefficients
$\partial^j Q_C : \bwedge^j_C \mfrak{g} \to \mfrak{g}$.
But giving such a homomorphism $\partial^j Q_C$ is the same as 
giving a $C$-multilinear homomorphism
$\partial^j Q : \bwedge^j \mfrak{g} \to \mfrak{g}$,
so there is a corresponding $C$-multilinear coderivation
$Q : \mrm{S}(\mfrak{g}[1]) \to \mrm{S}(\mfrak{g}[1])$.
One checks that $Q \circ Q = 0$ iff $Q_C \circ Q_C = 0$. 

Similarly for coalgebra homomorphisms. 
\end{proof}

An element 
$\gamma \in \mrm{S}_C(\mfrak{g}[1])$
is called {\em primitive} if
$\Delta_C(\gamma) = \gamma \otimes 1 + 1 \otimes \gamma$. 

\begin{lem} \label{lem9.2}
The set of primitive elements of $\mrm{S}_C(\mfrak{g}[1])$
is precisely 
$\mrm{S}_C^1(\mfrak{g}[1]) = \mfrak{g}[1]$. 
\end{lem}

\begin{proof}
By definition of the comultiplication any 
$\gamma \in \mfrak{g}[1]$
is primitive. For the converse, let us denote by $\mu$ 
the multiplication in $\mrm{S}_C(\mfrak{g}[1])$. 
One checks that 
$(\mu \circ \Delta_C)(\gamma) = 2^i \gamma$
for $\gamma \in \mrm{S}^i_C(\mfrak{g}[1])$. 
If $\gamma$ is primitive then 
$(\mu \circ \Delta_C)(\gamma) = 2 \gamma$,
so indeed $\gamma \in \mrm{S}^1_C(\mfrak{g}[1])$.
\end{proof}

Now let's assume that $C$ is a local ring, with nilpotent maximal 
ideal $\mfrak{m}$. Suppose we are given two $C$-multilinear 
$\mrm{L}_{\infty}$ algebras $(\mfrak{g}, Q)$ and 
$(\mfrak{g}', Q')$, and a $C$-multilinear $\mrm{L}_{\infty}$ 
morphism 
$\Psi : (\mfrak{g}, Q) \to (\mfrak{g}', Q')$.
Because the coderivation $Q$ is $C$-multilinear, 
the $C$-submodule
$\mfrak{m} \mfrak{g} \subset \mfrak{g}$
becomes a $C$-multilinear $\mrm{L}_{\infty}$ algebra
$(\mfrak{m} \mfrak{g}, Q)$. 
Likewise for $\mfrak{m} \mfrak{g}'$, and
$\Psi : (\mfrak{m} \mfrak{g}, Q) \to (\mfrak{m} \mfrak{g}', Q')$
is a $C$-multilinear $\mrm{L}_{\infty}$ morphism.

The fact that 
$\mfrak{m}$ is nilpotent is essential for the next definition.

\begin{dfn} \label{dfn9.1}
The {\em Maurer-Cartan equation} in 
$(\mfrak{m} \mfrak{g}, Q)$ is
\[ \sum_{i = 1}^{\infty} 
\smfrac{1}{i!} (\partial^i Q)(\omega^i) = 0 \]
for 
$\omega \in (\mfrak{m} \mfrak{g})^1 = 
(\mfrak{m} \mfrak{g}[1])^0$. 
\end{dfn}

An element $e \in \mrm{S}_C(\mfrak{g}[1])$ is called 
{\em group-like} if $\Delta_C(e) = e \otimes e$. For 
$\omega \in \mfrak{m} \mfrak{g}^1$
we define
\[ \opn{exp}(\omega) := \sum_{i \geq 0}
\smfrac{1}{i!} \omega^i \in \mrm{S}_C(\mfrak{g}[1]) . \]

\begin{lem} \label{lem9.3}
The function $\opn{exp}$ is a bijection from 
$\mfrak{m} \mfrak{g}[1]$
to the set of invertible group-like elements of 
$\mrm{S}_C(\mfrak{g}[1])$, with inverse $\opn{ln}$. 
\end{lem}

\begin{proof}
Let $\omega \in \mfrak{m} \mfrak{g}[1]$
and $e := \opn{exp}(\omega)$. The element $e$ is invertible, with 
inverse $\opn{exp}(-\omega)$. Using the fact that 
$\Delta_C(\omega) = \omega \otimes 1 + 1 \otimes \omega$
it easily follows that 
$\Delta_C(e) = e \otimes e$.
And trivially $\opn{ln}(e) = \omega$. 

For the opposite direction, let $e$ be invertible and group-like. 
Write it as $e = \sum_i \gamma_i$, with 
$\gamma_i \in \mrm{S}^i_C(\mfrak{g}[1])$.
Since $e$ is invertible one must have
$\gamma_0 \in C - \mfrak{m}$, and 
$\gamma_i \in \mfrak{m} \mrm{S}^i_C(\mfrak{g}[1])$
for all $i \geq 1$. The equation $\Delta_C(e) = e \otimes e$ 
implies that
\[ \Delta_C(\gamma_i) = \sum_{j +k = i} \gamma_j \otimes \gamma_k \]
for all $i$. Hence
\begin{equation} \label{eqn9.1}
2^i \gamma_i = \mu(\Delta_C(\gamma_i)) 
= \sum_{j +k = i} \gamma_j  \gamma_k .
\end{equation}
For $i = 0$ we get $\gamma_0 = \gamma_0^2$, and since $\gamma_0$
is invertible, it follows that $\gamma_0 = 1$. Let 
$\omega := \gamma_1 = \opn{ln}(e) \in 
\mfrak{m} \mrm{S}^1_C(\mfrak{g}[1]) =
\mfrak{m} \mfrak{g}[1]$. 
Using induction and equation (\ref{eqn9.1}) we see that 
$\gamma_i = \smfrac{1}{i!} \omega^i$ for all $i$. 
Thus $e = \opn{exp}(\omega)$. 
\end{proof}

\begin{lem} \label{lem9.4}
Let $\omega \in (\mfrak{m} \mfrak{g}[1])^0 = \mfrak{m} \mfrak{g}^1$ 
and
$e := \opn{exp}(\omega)$. Then $\omega$ is a solution of the MC 
equation iff $Q(e) = 0$. 
\end{lem}

\begin{proof}
Since $e$ is group-like  and invertible
(by Lemma \ref{lem9.3}) we have
\[ \Delta_C(Q(e)) = Q(e) \otimes e + e \otimes Q(e)  \]
and
\[ \Delta_C(e^{-1} Q(e)) = 
\Delta_C(e)^{-1} \Delta_C(Q(e)) =
e^{-1} Q(e) \otimes 1 + 1 \otimes e^{-1} Q(e) . \]
So the element $e^{-1} Q(e)$ is primitive, and by Lemma \ref{lem9.2} 
we get 
$e^{-1} Q(e) \in \mfrak{g}[1]$. 
On the other hand hence $Q(e)$ has no $0$-order term, and 
$Q(1) = 0$. Thus in the $1$st order term we get
\begin{equation} \label{eqn9.2}
\begin{aligned}
e^{-1} Q(e) & = \opn{ln} \bigl( e^{-1} Q(e) \bigr) \\
& = \opn{ln} \bigl( (1 - \omega + \smfrac{1}{2} \omega^2 \pm \cdots)
Q(e) \bigr) \\
& = \opn{ln} \bigl( Q(e) \bigr) \\
& = \sum_{i=0}^{\infty} 
\smfrac{1}{i!} \opn{ln} \bigl( Q(\omega^i) \bigr) \\
& = \sum_{i=1}^{\infty} \smfrac{1}{i!} (\partial^i Q)(\omega^i) . 
\end{aligned}
\end{equation}
Since $e$ is invertible we are done.
\end{proof}

\begin{lem} \label{lem9.5}
Given an element $\omega \in \mfrak{m} \mfrak{g}[1]$,
define
$\omega' := \sum_{i=1}^{\infty} \smfrac{1}{i!} 
(\partial^i \Psi)(\omega^i)
\in \mfrak{m} \mfrak{g}'[1]$,
$e := \opn{exp}(\omega)$ and
$e' := \opn{exp}(\omega')$. Then
$e' = \Psi(e)$. 
\end{lem}

\begin{proof}
 From Lemma \ref{lem9.3} we see that 
$\Delta_C(e) = e \otimes e$, and therefore also 
$\Delta_C(\Psi(e)) = \Psi(e) \otimes \Psi(e)
\in \mrm{S}_C(\mfrak{g}'[1])$. 
Since $\Psi$ is $C$-linear and $\Psi(1) = 1$ we get
$\Psi(e) \in 1 + \mfrak{m} \mrm{S}(\mfrak{g}'[1])$. 
Thus $\Psi(e)$ is group-like and invertible. 
According to Lemma \ref{lem9.3} it suffices to prove that 
$\opn{ln}(e') = \opn{ln}(\Psi(e))$. 
Now $\opn{ln}(e') = \omega'$ by definition. Since $\Psi(1)  = 1$ 
and $\opn{ln}(1) = 0$ it follows that
\[ \opn{ln}(\Psi(e)) = 
\opn{ln}(\Psi \big( \sum_{i=0}^{\infty} 
\smfrac{1}{i!} \omega^i \big) )
= \sum_{i=0}^{\infty} \smfrac{1}{i!} \opn{ln}(\Psi(\omega^i))
= \sum_{i=1}^{\infty} \smfrac{1}{i!} (\partial^i \Psi)(\omega^i)
= \omega' . \]
\end{proof}

\begin{prop} \label{prop9.3}
Suppose 
$\omega \in \mfrak{m} \mfrak{g}^1$
is a solution of the MC equation in
$(\mfrak{m} \mfrak{g}, Q)$. Define
$\omega' := \sum_{i=1}^{\infty} \smfrac{1}{i!} 
(\partial^i \Psi)(\omega^i)
\in \mfrak{m} \mfrak{g}'^1$.
Then $\omega'$ is a solution of the MC equation in 
$(\mfrak{m} \mfrak{g}', Q')$.
\end{prop}

\begin{proof}
Let $e := \opn{exp}(\omega)$ and
$e' := \opn{exp}(\omega')$. By Lemma \ref{lem9.4} we get
$Q(e) = 0$. Hence $Q'(\Psi(e)) = \Psi(Q(e)) = 0$. 
According to Lemma \ref{lem9.5} we have
$\Psi(e) = e'$, so $Q'(e') = 0$. Again by Lemma \ref{lem9.4} we 
deduce that $\omega'$ solves the MC equation. 
\end{proof}

\begin{dfn} \label{dfn9.2}
Let $\omega \in \mfrak{m} \mfrak{g}^1$. 
\begin{enumerate}
\item The coderivation $Q_{\omega}$ of 
$\mrm{S}_C(\mfrak{g}[1])$ over $C$, with $Q_{\omega}(1) := 0$
and with Taylor coefficients
\[ (\partial^i Q_{\omega})(\gamma) :=
\sum_{j \geq 0} \smfrac{1}{j!} (\partial^{i+j} Q)
(\omega^j \gamma)  \]
for $i \geq 1$ and $\gamma \in \mrm{S}^i_C(\mfrak{g}[1])$,
is called the {\em twist of $Q$ by $\omega$}.
\item The coalgebra homomorphism 
$\Psi_{\omega} : \mrm{S}_C(\mfrak{g}[1]) \to
\mrm{S}_C(\mfrak{g}'[1])$
over $C$, with $\Psi_{\omega}(1) := 1$ and
Taylor coefficients
\[ (\partial^i \Psi_{\omega})(\gamma) :=
\sum_{j \geq 0} \smfrac{1}{j!} (\partial^{i+j} \Psi)
(\omega^j \gamma)  \]
for $i \geq 1$ and $\gamma \in \mrm{S}_C^i(\mfrak{g}[1])$,
is called the {\em twist of $\Psi$ by $\omega$}.
\end{enumerate}
\end{dfn}

\begin{thm} \label{thm9.1}
Let $C$ be a commutative local $\K$-algebra with nilpotent maximal 
ideal $\mfrak{m}$. Let $(\mfrak{g}, Q)$ and $(\mfrak{g}', Q')$ 
be $C$-multilinear $\mrm{L}_{\infty}$ 
algebras and $\Psi : (\mfrak{g}, Q) \to (\mfrak{g}', Q')$ 
a $C$-multilinear $\mrm{L}_{\infty}$ morphism.
Suppose $\omega \in \mfrak{m} \mfrak{g}^1$ 
a solution of the MC equation in 
$(\mfrak{m} \mfrak{g}, Q)$. Define
\[ \omega' := \sum_{i=1}^{\infty} \smfrac{1}{j!}
(\partial^j \Psi)(\omega^j) \in \mfrak{m} {\mfrak{g}'}^1 . \] 
Then $(\mfrak{g}, Q_{\omega})$ and
$(\mfrak{g}', Q'_{\omega'})$ 
are $\mrm{L}_{\infty}$ algebras, and
\[ \Psi_{\omega} : (\mfrak{g}, Q_{\omega}) \to
(\mfrak{g}', Q'_{\omega'}) \]
is an $\mrm{L}_{\infty}$ morphism.
\end{thm}

\begin{proof}
Let $e := \opn{exp}(\omega)$. Define 
$\Phi_e : \mrm{S}_C(\mfrak{g}[1]) \to \mrm{S}_C(\mfrak{g}[1])$
to be $\Phi_e(\gamma) := e \gamma$. Since $e$ is group-like and
invertible it 
follows that $\Phi_e$ is a coalgebra automorphism. Therefore
$\til{Q}_{\omega} := \Phi_e^{-1} \circ Q \circ \Phi_e$
is a degree $1$ coderivation of $\mrm{S}_C(\mfrak{g}[1])$,
satisfying 
$\til{Q}_{\omega} \circ \til{Q}_{\omega} = 0$
and
$\til{Q}_{\omega}(1) = e^{-1} Q(e) = 0$;
cf.\ Lemma \ref{lem9.4}. So 
$(\mfrak{g}, \til{Q}_{\omega})$
is an $\mrm{L}_{\infty}$ algebra. Likewise we have a 
coalgebra automorphism $\Phi_{e'}$ and a coderivation
$\til{Q}'_{\omega'} := \Phi_{e'}^{-1} \circ Q'
\circ \Phi_{e'}$
of $\mrm{S}_C(\mfrak{g}'[1])$,
where $e' := \opn{exp}(\omega')$.
The degree $0$ coalgebra homomorphism
$\til{\Psi}_{\omega} := \Phi_{e'}^{-1} \circ \Psi \circ \Phi_{e}$
satisfies
$\til{\Psi}_{\omega} \circ \til{Q}_{\omega} = 
\til{Q}'_{\omega'} \circ \til{\Psi}_{\omega}$, 
and also 
$\til{\Psi}_{\omega}(1) = {e'}^{-1} \Psi(e) = 
{e'}^{-1} e' = 1$,
by Lemma \ref{lem9.5}. Hence we have an $\mrm{L}_{\infty}$ morphism
$\til{\Psi}_{\omega} : 
(\mfrak{g}, \til{Q}_{\omega}) \to (\mfrak{g}', \til{Q}'_{\omega'})$.

Let us calculate the Taylor coefficients of 
$\til{Q}_{\omega}$. 
For $\gamma \in \mrm{S}_C^i(\mfrak{g}[1])$ one has
\[ (\partial^i \til{Q}_{\omega})(\gamma) = 
\opn{ln}(\til{Q}_{\omega}(\gamma)) = 
\opn{ln}(e^{-1} Q (e \gamma)) . \]
But just as in (\ref{eqn9.2}), since $Q (e \gamma)$
has no zero order term, we obtain
\[ \opn{ln}(e^{-1} Q (e \gamma)) = \opn{ln}(Q (e \gamma)) . \]
And
\begin{equation} \label{eqn9.3}
\begin{aligned}
\opn{ln}(Q (e \gamma)) 
& = \opn{ln} \bigl( Q \Bigl( \sum\nolimits_{j \geq 0} 
\smfrac{1}{j!} \omega^j \gamma \Bigr) \bigr)  \\
& = \sum_{j \geq 0} \smfrac{1}{j!}
\opn{ln}(Q(\omega^j \gamma)) \\
& = \sum_{j \geq 0} \smfrac{1}{j!}
(\partial^{i+j} Q)(\omega^j \gamma) \\
& = (\partial^i Q_{\omega})(\gamma) .
\end{aligned} 
\end{equation}
Therefore $\til{Q}_{\omega} = Q_{\omega}$. 
Similarly we see that 
$\til{Q}'_{\omega'} = Q'_{\omega'}$ 
and
$\til{\Psi}_{\omega} = \Psi_{\omega}$. 
\end{proof}

\begin{rem}
The formulation of Theorem \ref{thm9.1}, 
as well as the idea for the proof, were suggested by 
Vladimir Hinich. An analogous result, for $\mrm{A}_{\infty}$ 
algebras, is in \cite[Section 6.1]{Le}.
\end{rem}

If $(\mfrak{g}, Q)$ is a DG Lie algebra then the sum occurring in 
Definition \ref{dfn9.2}(1) is finite, so the coderivation 
$Q_{\omega}$ can be defined without a nilpotence assumption
on the coefficients.

\begin{lem} \label{lem9.6} 
Let $(\mfrak{g}, Q)$ be a DG Lie algebra, and let
$\omega \in \mfrak{g}^1$ be a solution of the MC equation. Then 
the $\mrm{L}_{\infty}$ algebra $(\mfrak{g}, Q_{\omega})$ is also a 
DG Lie algebra. In fact, for $\gamma_i \in \mfrak{g}$ one has
\[ (\partial^1 Q_{\omega})(\gamma_1) = 
(\partial^1 Q)(\gamma_1) + (\partial^2 Q)(\omega \gamma_1) =
\mrm{d}(\gamma_1) + [\omega, \gamma_1] =
(\d + \opn{ad}(\omega))(\gamma_1) , \]
\[ (\partial^2 Q_{\omega})(\gamma_1 \gamma_2) = 
(\partial^2 Q)(\gamma_1 \gamma_2) = [\gamma_1, \gamma_2] , \]
and $\partial^j Q_{\omega} = 0$ for $j \geq 3$.
\end{lem}

\begin{proof}
Like equation (\ref{eqn9.3}), with $C := \K$ and $e := 1$.
\end{proof}

In the situation of the lemma, the twisted DG Lie algebra 
$(\mfrak{g}, Q_{\omega})$ will usually be denoted by
$\mfrak{g}_{\omega}$.

Let $A$ be a DG super-commutative associative unital DG algebra
in \linebreak $\cat{Dir} \cat{Inv} \cat{Mod} \K$.
The notion of DG $A$-module Lie algebra 
in $\cat{Dir} \cat{Inv} \cat{Mod} \K$
was introduced in Definition \ref{dfn1.4}.

\begin{dfn}
Let $A$ be a DG super-commutative associative unital 
DG algebra in $\cat{Dir} \cat{Inv} \cat{Mod} \K$,
let $\mfrak{g}$ and $\mfrak{g}'$ be DG $A$-module 
Lie algebras in $\cat{Dir} \cat{Inv} \cat{Mod} \K$, and let
$\Psi : \mfrak{g} \to \mfrak{g}'$ be an $\mrm{L}_{\infty}$ morphism.
\begin{enumerate}
\item If each Taylor coefficient 
$\partial^j \Psi : \prod^j \mfrak{g} \to \mfrak{g}'$ 
is continuous then we say that $\Psi$ is a {\em continuous 
$\mrm{L}_{\infty}$ morphism}. 
\item Assume each Taylor coefficient 
$\partial^j \Psi : \prod^j \mfrak{g} \to \mfrak{g}'$ 
is $A$-multilinear, i.e.\ 
\[ (\partial^j \Psi)(a_1 \gamma_1, \ldots, a_j \gamma_j) = 
\pm a_1 \cdots a_j \cdot 
(\partial^j \Psi)(\gamma_1, \ldots, \gamma_j) \]
for all homogeneous elements $a_k \in A$ and 
$\gamma_k \in \mfrak{g}$,
with sign according to the Koszul rule, then we say that 
$\Psi$ is an {\em $A$-multilinear $\mrm{L}_{\infty}$ morphism}. 
\end{enumerate}
\end{dfn}

\begin{prop} \label{prop10.1} 
Let $A$ and $B$ be DG super-comm\-utative associative unital 
DG algebras in $\cat{Dir} \cat{Inv} \cat{Mod} \K$,
and let $\mfrak{g}$ and $\mfrak{g}'$ 
be DG $A$-module Lie algebras in 
$\cat{Dir} \cat{Inv} \cat{Mod} \K$.
Suppose $A \to B$ is a continuous DG algebra homomorphism, and 
$\Psi : \mfrak{g} \to \mfrak{g}'$ 
is a continuous $A$-multilinear $\mrm{L}_{\infty}$ morphism.
Let 
$\partial^j \Psi_{\what{B}} : \prod^j 
(B \hatotimes{A} \mfrak{g}) \to 
B \hatotimes{A} \mfrak{g}'$
be the unique continuous $\what{B}$-multilinear homomorphism
extending $\partial^j \Psi$.
Then the degree $0$ coalgebra homomorphism 
\[ \Psi_{\what{B}} : \mrm{S}(B \hatotimes{A} \mfrak{g}[1]) \to 
\mrm{S}(B \hatotimes{A} \mfrak{g}'[1]) , \]
with $\Psi_{\what{B}}(1) := 1$ and with Taylor coefficients 
$\partial^j \Psi_{\what{B}}$, is an $\mrm{L}_{\infty}$ morphism
\[ \Psi_{\what{B}} : B \hatotimes{A} \mfrak{g} \to 
B \hatotimes{A} \mfrak{g}' . \]
\end{prop}

\begin{proof}
First consider the continuous $B$-multilinear homomorphisms
$\partial^j \Psi_{B} :$ \linebreak $\prod^j 
(B \otimes_{A} \mfrak{g}) \to B \otimes_{A} \mfrak{g}'$
extending $\partial^j \Psi$.
It is a straightforward calculation to verify that the 
$\mrm{L}_{\infty}$ morphism identities of Proposition 
\ref{prop9.1} hold for the sequence of operators
$\set{ \partial^j \Psi_{B} }_{j \geq 1}$.
The completion process respects these identities (cf.\ proof of 
Proposition \ref{prop1.16}).
\end{proof}

\begin{thm} \label{thm9.2}
Let $\mfrak{g}$ and $\mfrak{g}'$ be DG Lie algebras
in $\cat{Dir} \cat{Inv} \cat{Mod} \K$, and let
$\Psi : \mfrak{g} \to \mfrak{g}'$ 
be a continuous $\mrm{L}_{\infty}$ morphism. 
Let $A = \boplus_{i \in \mbb{N}} A^i$ be a complete
associative unital super-commutative DG
algebra in $\cat{Dir} \cat{Inv} \cat{Mod} \K$. By Proposition 
\tup{\ref{prop10.1}} there is an induced continuous $A$-multilinear 
$\mrm{L}_{\infty}$ morphism
$\Psi_{A} : A \hatotimes{} \mfrak{g} 
\to A \hatotimes{} \mfrak{g}'$.
Let $\omega \in A^1 \hatotimes{} \mfrak{g}^0$ 
be a solution of the MC equation in $A \hatotimes{} \mfrak{g}$.
Assume $\d_{\mfrak{g}} = 0$, 
$(\partial^j \Psi_A)(\omega^j) = 0$ for all $j \geq 2$,
and also that $\mfrak{g}'$ is bounded below. Define 
$\omega' := (\partial^1 \Psi_A)(\omega) \in 
A^1 \hatotimes{} {\mfrak{g}'}^0$.
Then:
\begin{enumerate}
\item The element $\omega'$ is a solution of the MC equation in
$A \hatotimes{} \mfrak{g}'$.
\item Given 
$c \in \mrm{S}^j \big( A \hatotimes{} \mfrak{g}[1] \big)$
there exists a natural number $k_0$ such that \linebreak
$(\partial^{j+k} \Psi_{A})(\omega^k c) = 0$
for all $k > k_0$. 
\item The degree $0$ coalgebra homomorphism 
\[ \Psi_{A, \omega} :
\mrm{S} \big( A \hatotimes{} \mfrak{g}[1] \big)
\to
\mrm{S} \big( A \hatotimes{} \mfrak{g}'[1] \big) ,
\]
with 
$\Psi_{A, \omega}(1) := 1$
and Taylor coefficients
\[ (\partial^j \Psi_{A, \omega}) (c) :=
\sum_{k \geq 0} \smfrac{1}{(j+k)!}
(\partial^{j+k} \Psi_{A}) \big( \omega^k c \big) \]
for 
$c \in \mrm{S}^j \big( A \hatotimes{} \mfrak{g}[1] \big)$,
is a continuous $A$-multilinear $\mrm{L}_{\infty}$ morphism
\[ \Psi_{A, \omega} :
\big( A \hatotimes{} \mfrak{g} \big)_{\omega} \to
\big( A \hatotimes{} \mfrak{g}' \big)_{\omega'} . \]
\end{enumerate}
\end{thm}

\begin{proof}
We shall use a ``deformation argument''. Consider the base field 
$\K$ as a discrete inv $\K$-module. The polynomial 
algebra $\K[\hbar]$ is endowed with the dir-inv $\K$-module 
structure such that the homomorphism
$\boplus_{i \in \mbb{N}} \K \to \K[\hbar]$,
whose $i$-th component is multiplication by $\hbar^i$,
is an isomorphism in $\cat{Dir} \cat{Inv} \cat{Mod} \K$.
Note that $\K[\hbar]$ is a discrete dir-inv module, but it is not 
trivial. We view $\K[\hbar]$ as a DG algebra concentrated in 
degree $0$ (with zero differential).

For any $i \in \mbb{N}$ let  
$A[\hbar]^i := \K[\hbar] \otimes{} A^i$,
and let
$A[\hbar] := \boplus_{i \in \mbb{N}} A[\hbar]^i$,
which is a DG algebra 
in $\cat{Dir} \cat{Inv} \cat{Mod} \K$, with differential 
$\d_{A[\hbar]} := \bsym{1} \otimes \d_A$.
We will need a ``twisted'' version of $A[\hbar]$, which we denote 
by $A[\hbar]^{\sim}$. Let
$A[\hbar]^{\sim \, i} := \hbar^i A[\hbar]^i$,
and define 
$A[\hbar]^{\sim} := \boplus_{i \in \mbb{N}}  A[\hbar]^{\sim \, i}$,
which his a graded subalgebra of $A[\hbar]$. The differential is 
$\d_{A[\hbar]^{\sim}} := \hbar \d_{A[\hbar]}$.
The dir-inv structure is such that the homomorphism
$\boplus_{i, j \in \mbb{N}} A^i \to A[\hbar]^{\sim}$,
whose $(i, j)$-th component is multiplication by $\hbar^{i+j}$,
is an isomorphism in $\cat{Dir} \cat{Inv} \cat{Mod} \K$.
The specialization $\hbar \mapsto 1$ is a continuous DG algebra 
homomorphism $A[\hbar]^{\sim} \to A$.
There is an induced continuous $A[\hbar]^{\sim}$-multilinear 
$\mrm{L}_{\infty}$ morphism
$\Psi_{A[\hbar]^{\sim}} : A[\hbar]^{\sim} \hatotimes{} \mfrak{g} 
\to A[\hbar]^{\sim} \hatotimes{} \mfrak{g}'$.

We proceed in several steps.

\medskip \noindent
Step 1. Say $r_0$ bounds $\mfrak{g}'$ from below, i.e.\ 
${\mfrak{g}'}^r = 0$ for all $r < r_0$. Take some $j \geq 1$. 
For any $l \in \set{1, \ldots, j}$ choose 
$p_l, q_l \in \mbb{Z}$, $\gamma_l \in \mfrak{g}^{p_l}$ and
$a_l \in A[\hbar]^{\sim\, q_l}$. 
Also choose $\gamma_0 \in \mfrak{g}^{0}$ and
$a_0 \in A[\hbar]^{\sim\, 1}$. Let $p := \sum_{l=1}^j p_l$ and
$q := \sum_{l=1}^j q_l$. Because 
$\partial^{j+k} \Psi_{A[\hbar]^{\sim}}$ 
is induced from $\partial^{j+k} \Psi$, and this is a homogeneous map 
of degree $1-j-k$, we have
\[ \begin{aligned}
& (\partial^{j+k} \Psi_{A[\hbar]^{\sim}})
\big( (a_0 \otimes \gamma_0)^k (a_1 \otimes \gamma_1)
\cdots (a_j \otimes \gamma_j) \big) \\
& \qquad = \pm a_0^k a_1 \cdots a_j \otimes 
(\partial^{j+k} \Psi)
(\gamma_0^k \gamma_1 \cdots \gamma_j) 
\in A[\hbar]^{\sim \, k+q} \hatotimes{} \mfrak{g}^{p+1-j-k} . 
\end{aligned} \]
But $\mfrak{g}^{p+1-j-k} = 0$ for all $k > p+1-j-r_0$.

Using multilinearity and continuity we conclude that given any
$c \in$ \linebreak 
$\mrm{S}^j \big( A[\hbar]^{\sim} \hatotimes{} \mfrak{g}[1] \big)$
there exists a natural number $k_0$ such that 
$(\partial^{j+k} \Psi_{A[\hbar]^{\sim}}) 
\big( (\hbar \omega)^k c \big)$ \linebreak $= 0$
for all $k > k_0$.

\medskip \noindent
Step 2. We are going to prove that $\hbar \omega$ 
is a solution of the MC equation in 
$A[\hbar]^{\sim} \hatotimes{} \mfrak{g}$.
It is given that $\omega$ is a solution of the MC equation in 
$A \hatotimes{} \mfrak{g}$. Because $\d_{\mfrak{g}} = 0$, this 
means that 
\[ (\mrm{d}_{A} \otimes \bsym{1})(\omega) +
\smfrac{1}{2} [\omega, \omega] = 0 . \]
Hence
\[ \mrm{d}_{A[\hbar]^{\sim} \hatotimes{} \mfrak{g}}(\hbar \omega) +
\smfrac{1}{2} [\hbar \omega, \hbar \omega] =
\hbar^2 (\mrm{d}_{A} \otimes \bsym{1})(\omega) +
\smfrac{1}{2} \hbar^2 [\omega, \omega] = 0 . \]
So $\hbar \omega$ solves the MC equation in 
$A[\hbar]^{\sim} \hatotimes{} \mfrak{g}$. 

\medskip \noindent
Step 3. Now we shall prove that $\hbar \omega'$ solves the MC 
equation in $A[\hbar]^{\sim} \hatotimes{} \mfrak{g}'$.
This will require an infinitesimal argument. For any natural number 
$m$ define $\K[\hbar]_m := \K[\hbar] / (\hbar^{m+1})$ and
$A [\hbar]_m := \K[\hbar]_m \otimes{} A$.
The latter is a DG algebra with differential
$\d_{A [\hbar]_m} := \bsym{1} \otimes \d_A$.
Let 
$A [\hbar]_m^{\sim} := \boplus_{i=0}^m \hbar^i  A[\hbar]_m^i$,
which is a subalgebra of $A[\hbar]_m$, but its differential is
$\d_{A [\hbar]_m^{\sim}} := \hbar \d_{A [\hbar]_m}$.
There is a surjective DG Lie algebra homomorphism
$A[\hbar]^{\sim} \hatotimes{} \mfrak{g}' \to 
A[\hbar]^{\sim}_m \hatotimes{} \mfrak{g}'$,
with kernel
$\big( A[\hbar]^{\sim} \cap \hbar^{m+1} A [\hbar] \big)
\hatotimes{} \mfrak{g}'$.
Since $\bigcap_{m \geq 0} \hbar^{m+1} A[\hbar] = 0$,
it it suffices to prove that $\hbar \omega'$ solves the 
MC equation in $A[\hbar]_m^{\sim} \hatotimes{} \mfrak{g}'$.

Now $C := \K[\hbar]_m$ is an artinian local ring with maximal ideal 
$\mfrak{m} := (\hbar)$. Define the DG Lie algebra 
$\mfrak{h} := A[\hbar]_m \hatotimes{} \mfrak{g}$, with differential
$\d_{\mfrak{h}} :=
\hbar \d_{A[\hbar]_m} \otimes \bsym{1} + \bsym{1} \otimes
\d_{\mfrak{g}}$;
so $A[\hbar]_m^{\sim} \hatotimes{} \mfrak{g} \subset \mfrak{h}$
as DG Lie algebras. Similarly define $\mfrak{h}'$. 
There is a $C$-multilinear $\mrm{L}_{\infty}$ morphism 
$\Phi : \mfrak{h} \to \mfrak{h}'$ extending
$\Psi_{A[\hbar]_m^{\sim}} : A[\hbar]^{\sim}_m \hatotimes{} \mfrak{g} 
\to A[\hbar]^{\sim}_m \hatotimes{} \mfrak{g}'$.
By step 2 the element $\nu := \hbar \omega \in \mfrak{m} \mfrak{h}$
is a solution of the MC equation. According to 
Proposition \ref{prop9.3}
the element
$\nu' := \sum_{k \geq 1} (\partial^{k} \Phi)(\nu^k)$
is a solution of the MC equation in $\mfrak{h}'$.
But $\nu' = \hbar \omega'$.

\medskip \noindent
Step 4. Pick a natural number $m$. Let 
$\mfrak{h}, \mfrak{h}', \Phi, \nu$ and $\nu'$ be as in step 3. 
According to Theorem \ref{thm9.1} there is a twisted 
$\mrm{L}_{\infty}$ morphism
$\Phi_{\nu} : \mfrak{h}_{\nu} \to \mfrak{h}'_{\nu'}$.
Since 
$(A[\hbar]^{\sim}_m \hatotimes{} \mfrak{g})_{\hbar \omega} \subset 
\mfrak{h}_{\nu}$
and
$(A[\hbar]^{\sim}_m \hatotimes{} \mfrak{g}')_{\hbar \omega'} \subset 
\mfrak{h}'_{\nu'}$
as DG Lie algebras, and $\Phi_{\nu}$ extends 
$\Psi_{A[\hbar]^{\sim}_m, \hbar \omega}$, 
it follows that
$\Psi_{A[\hbar]^{\sim}_m, \hbar \omega} : 
A[\hbar]^{\sim}_m \hatotimes{} \mfrak{g} \to 
A[\hbar]^{\sim}_m \hatotimes{} \mfrak{g}'$
is an $\mrm{L}_{\infty}$ morphism. 
This means that the Taylor coefficients
\[ \partial^j \Psi_{A[\hbar]^{\sim}_m, \hbar \omega} : 
\prod\nolimits^j 
(A[\hbar]^{\sim}_m \hatotimes{} \mfrak{g})_{\hbar \omega} \to 
(A[\hbar]^{\sim}_m \hatotimes{} \mfrak{g}')_{\hbar \omega'} \]
satisfy the identities of Proposition \ref{prop9.1}.
As explained in step 3, this implies that 
\[ \partial^j \Psi_{A[\hbar]^{\sim}, \hbar \omega} : 
\prod\nolimits^j 
(A[\hbar]^{\sim} \hatotimes{} \mfrak{g})_{\hbar \omega} \to 
(A[\hbar]^{\sim} \hatotimes{} \mfrak{g}')_{\hbar \omega'} \]
also satisfy these identities. We conclude that
$\Psi_{A[\hbar]^{\sim}, \hbar \omega}$ is an $\mrm{L}_{\infty}$ 
morphism.  

\medskip \noindent
Step 5. Specialization $\hbar \mapsto 1$ induces surjective DG 
Lie algebra homomorphisms
$A[\hbar]^{\sim} \hatotimes{} \mfrak{g} \to 
A \hatotimes{} \mfrak{g}$
and
$A[\hbar]^{\sim} \hatotimes{} \mfrak{g}' 
\to A \hatotimes{} \mfrak{g}'$,
sending $\hbar \omega \mapsto \omega$,
$\hbar \omega' \mapsto \omega'$ and
$\Psi_{A[\hbar]^{\sim}, \hbar \omega} \mapsto
\Psi_{A, \omega}$.
Therefore assertions (1-3) of the theorem hold.
\end{proof}

\section{The Universal $\mrm{L}_{\infty}$ Morphism
of Kontsevich}
\label{sec4}

In this section $\K$ is a field of characteristic $0$ and 
$C$ is a commutative $\K$-algebra. Recall that we denote by 
$\mcal{T}_{C} = \mcal{T}(C / \K) := \opn{Der}_{\K}(C)$,
the module of derivations of $C$ relative to 
$\K$. This is a Lie algebra over $\K$. 
Following \cite{Ko1} we make the next definitions.

\begin{dfn} \label{dfn1.1}
For $p \geq -1$ let 
\[ \mcal{T}^{p}_{\mrm{poly}}(C) :=
\bwedge^{p+1}_{C} \mcal{T}_{C} , \]
the module of {\em poly derivations} (or {\em poly tangents}) 
of degree $p$ of $C$ relative to $\K$. Let
\[ \mcal{T}_{\mrm{poly}}(C)  := 
\boplus_{p}\, \mcal{T}^{p}_{\mrm{poly}}(C). \] 
This is a DG Lie algebra, with zero differential, and with 
the Schouten-Nijenhuis bracket, which is determined by the 
formulas 
\[ \begin{aligned} 
{} [\alpha_1 \wedge \alpha_2, \alpha_3]  
& = \alpha_1 \wedge [\alpha_2, \alpha_3] + (-1)^{(p_2 + 1) p_3}
[\alpha_1, \alpha_3] \wedge \alpha_2 \\
& \text{and} \\
[\alpha_1 , \alpha_2] & = 
(-1)^{1 + p_1 p_2} [\alpha_2 , \alpha_1] 
\end{aligned} \]
for elements
$\alpha_i \in \mcal{T}^{p_i}_{\mrm{poly}}(C)$.
\end{dfn}

\begin{dfn} \label{dfn2.2}
For any $p \geq -1$ let
$\mcal{D}^{p}_{\mrm{poly}}(C)$ 
be the set of $\K$-multilinear multi differential operators 
$\phi : C^{p+1} \to C$ (see Definition \ref{dfn1.5}). The direct sum 
\[ \mcal{D}_{\mrm{poly}}(C) :=
\boplus_p \, \mcal{D}^{p}_{\mrm{poly}}(C) \]
is a DG Lie algebra. The differential $\mrm{d}_{\mcal{D}}$
is the shifted Hochschild differential, and the Lie bracket is the
Gerstenhaber bracket (see \cite[Section 3.4.2]{Ko1}).
The elements of $\mcal{D}_{\mrm{poly}}(C)$ 
are called {\em poly differential operators} 
relative to $\K$. 
\end{dfn}

In the notation of Section \ref{sec2} and Example \ref{exa1.2}
one has
\[ \mcal{D}^p_{\mrm{poly}}(C) = 
\mcal{D}\mathit{iff}_{\mrm{poly}}(C; 
\underset{p+1}{\underbrace{C, \ldots, C}}; C)
= \mcal{C}^{p+1}_{\mrm{cd}}(C) ; \]
see formula (\ref{eqn2.1}).

Observe that 
$\mcal{D}^{p}_{\mrm{poly}}(C) \subset
\opn{Hom}_{\K}(C^{\otimes (p+1)},  C)$,
and $\mcal{D}_{\mrm{poly}}(C)$ is a sub DG 
Lie algebra of the shifted Hochschild cochain complex of $C$ 
relative to $\K$. For $p = -1, 0$ we have
$\mcal{D}^{-1}_{\mrm{poly}}(C) = C$ and 
$\mcal{D}^{0}_{\mrm{poly}}(C) = \mcal{D}(C)$,
the ring of differential operators. Note that 
$\mcal{D}^{p}_{\mrm{poly}}(C)$ is a left module over 
$\mcal{D}(C)$, by the formula 
$D \cdot \phi := D \circ \phi$; 
and in this way it is also a left $C$-module.

When $C := \K[\bsym{t}] = \K[t_1, \ldots, t_n]$, the polynomial 
algebra in $n \geq 1$ variables, and $p \geq 1$,
the following is true. The $\K[\bsym{t}]$-module 
$\mcal{T}^{p-1}_{\mrm{poly}}(\K[\bsym{t}])$
is free with finite basis 
$\{ \smfrac{\partial}{\partial t_{i_1}} \wedge \cdots \wedge
\smfrac{\partial}{\partial t_{i_p}} \}$,
indexed by the sequences
$0 \leq i_1 < \cdots < i_p \leq n$. The $\K[\bsym{t}]$-module 
$\mcal{D}^{p-1}_{\mrm{poly}}(\K[\bsym{t}])$ is also free, with 
countable basis 
\begin{equation} \label{eqn2.9}
\{ (\smfrac{\partial}{\partial \bsym{t}})^{\bsym{j}_1} \otimes
\cdots \otimes 
(\smfrac{\partial}{\partial \bsym{t}})^{\bsym{j}_{p}} \}
_{\bsym{j}_1, \ldots, \bsym{j}_p \in \mbb{N}^n},  
\end{equation}
where for $\bsym{j}_k = (j_{k, 1}, \ldots, j_{k, n}) \in \mbb{N}^n$ 
we write
$(\frac{\partial}{\partial \bsym{t}})^{\bsym{j}_k} :=
(\frac{\partial}{\partial t_1})^{j_{k, 1}}
\cdots (\frac{\partial}{\partial t_n})^{j_{k, n}}$.

For any $p \geq -1$ let 
$\mrm{F}_m \mcal{D}^{p}_{\mrm{poly}}(C)$
be the set of poly differential operators of order $\leq m$ in 
each argument. This is $C$-submodule of 
$\mcal{D}^{p}_{\mrm{poly}}(C)$.

\begin{lem} \label{lem2.9}
\begin{enumerate}
\item For any $m, p$ one has
\[ \mrm{d}_{\mcal{D}} \big( 
\mrm{F}_m \mcal{D}^{p}_{\mrm{poly}}(C) \big) 
\subset
\mrm{F}_m \mcal{D}^{p+1}_{\mrm{poly}}(C) . \]
\item For any $m, m', p$ one has
\[ \big[ \mrm{F}_m \mcal{D}^{p}_{\mrm{poly}}(C) , 
\mrm{F}_{m'} \mcal{D}^{p'}_{\mrm{poly}}(C) \big] \subset 
\mrm{F}_{m+m'} \mcal{D}^{p+p'}_{\mrm{poly}}(C) ; \]
and
\[ [-,-] : \mrm{F}_m \mcal{D}^{p}_{\mrm{poly}}(C) \times
\mrm{F}_{m'} \mcal{D}^{p'}_{\mrm{poly}}(C) \to 
\mcal{D}^{p+p'}_{\mrm{poly}}(C) \]
is a poly differential operator of order $\leq m+m'$ in each of 
its two arguments.
\end{enumerate}
\end{lem}

\begin{proof}
These assertions follow easily
from the definitions of the Hochschild differential and the 
Gerstenhaber bracket; cf.\ \cite[Section 3.4.2]{Ko1}. 
\end{proof}

\begin{lem}
Assume $C$ is a finitely generated $\K$-algebra. 
Then $\mcal{T}^{p}_{\mrm{poly}}(C)$ and
$\mrm{F}_m \mcal{D}^{p}_{\mrm{poly}}(C)$
are finitely generated $C$-modules. 
\end{lem}

\begin{proof}
One has
\[ \mcal{T}^{p}_{\mrm{poly}}(C) \cong 
\opn{Hom}_A(\Omega^{p+1}_{C}, A) \] 
and
\[ \mrm{F}_m \mcal{D}^{p}_{\mrm{poly}}(C) \cong
\opn{Hom}_C \big( \mcal{C}_{p+1, m}(C), C \big) ; \]
see Lemma \ref{lem3.6}. 
The $C$-modules $\Omega^{p+1}_{C}$ and
$\mcal{C}_{p+1, m}(C)$ are finitely generated.
\end{proof}

\begin{prop} \label{prop3.8}
Assume $C$ is a finitely generated $\K$-algebra,
and $C'$ is a noetherian, $\mfrak{c}'$-adically complete, flat,
$\mfrak{c}'$-adically formally \'etale $C$-algebra. Let's write 
$\mcal{G}$ for either $\mcal{T}^{}_{\mrm{poly}}$ or
$\mcal{D}^{}_{\mrm{poly}}$. 
Then:
\begin{enumerate}
\item There is a DG Lie algebra homomorphism
$\mcal{G}(C) \to \mcal{G}(C')$, 
which is functorial in $C \to C'$.
\item The induced $C'$-linear homomorphism 
$C' \otimes_C \mcal{G}^{p}(C) \to \mcal{G}^{p}(C')$
is bijective.
\item For any $m$ the isomorphisms in \tup{(2)}, for
$\mcal{G} = \mcal{D}^{}_{\mrm{poly}}$, 
restrict to isomorphisms 
\[ C' \otimes_C \mrm{F}_m \mcal{D}^{p}_{\mrm{poly}}(C) \iso
\mrm{F}_m \mcal{D}^{p}_{\mrm{poly}}(C') . \]
\end{enumerate}
\end{prop}

\begin{proof}
Consider $\mcal{G} = \mcal{D}^{}_{\mrm{poly}}$.
Let $\phi \in \mcal{D}^{p}_{\mrm{poly}}(C)$. 
According to Proposition \ref{prop3.9}, applied to the case
$M_1, \ldots, M_{p+1}, N := A$, 
there is a unique 
$\phi' \in \mcal{D}^{p}_{\mrm{poly}}(C')$
extending $\phi$. From the definitions of the Gerstenhaber 
bracket and the Hochschild differential, it immediately follows 
that the function
$\mcal{D}^{}_{\mrm{poly}}(C) \to
\mcal{D}^{}_{\mrm{poly}}(C')$,
$\phi \mapsto \phi'$, 
is a DG Lie algebra homomorphism.
Parts (2,3) are also consequences of Proposition \ref{prop3.9}.

The case $\mcal{G} = \mcal{T}^{}_{\mrm{poly}}$ is
done similarly (and is well-known). 
\end{proof}

Consider $C := \K[\bsym{t}]$ and 
$C' := \K[[\bsym{t}]] = \K[[t_1, \ldots, t_n]]$, 
the power series algebra. Since 
$\mcal{T}^{p}_{\mrm{poly}}(\K[[\bsym{t}]]) \cong
\K[[\bsym{t}]] \otimes_{\K[\bsym{t}]}
\mcal{T}^{p}_{\mrm{poly}}(\K[\bsym{t}])$
is a finitely generated left $\K[[\bsym{t}]]$-module, it is an 
inv $\K[[\bsym{t}]]$-module with the $(\bsym{t})$-adic inv 
structure; cf.\ Example \ref{exa2.2}.
Likewise $\mcal{D}^{p}_{\mrm{poly}}(\K[[\bsym{t}]])$
is a dir-inv $\K[[\bsym{t}]]$-module. 
By Proposition \ref{prop3.8}, 
\[ \mrm{F}_m \mcal{D}^{p}_{\mrm{poly}}(\K[[\bsym{t}]]) \cong
\K[[\bsym{t}]] \otimes_{\K[\bsym{t}]}
\mrm{F}_m \mcal{D}^{p}_{\mrm{poly}}(\K[\bsym{t}]) , \]
which is a finitely generated $\K[[\bsym{t}]]$-module. So
according to Example \ref{exa1.5} we may take 
$\set{ \mrm{F}_m \mcal{D}^{p}_{\mrm{poly}}(\K[[\bsym{t}]]) }
_{m \in \mbb{N}}$
as the dir-inv structure of 
$\mcal{D}^{p}_{\mrm{poly}}(\K[[\bsym{t}]])$.
Now forgetting the $\K[[\bsym{t}]]$-module structure, 
$\mcal{T}^{p}_{\mrm{poly}}(\K[[\bsym{t}]])$
becomes an inv $\K$-module, and 
$\mcal{D}^{p}_{\mrm{poly}}(\K[[\bsym{t}]])$
becomes a dir-inv $\K$-module.

\begin{prop} \label{prop3.7} 
Let $\mcal{G}$ stand either for $\mcal{T}_{\mrm{poly}}$ or
$\mcal{D}_{\mrm{poly}}$. Then
$\mcal{G}(\K[[\bsym{t}]])$ 
is a complete DG Lie algebra in 
$\cat{Dir} \cat{Inv} \cat{Mod} \K$. 
\end{prop}

\begin{proof}
Use Proposition \ref{prop1.9}, and, for the case 
$\mcal{G} = \mcal{D}_{\mrm{poly}}$, also
Lemma \ref{lem2.9}.
\end{proof}

\begin{rem}
One might prefer to view 
$\mcal{T}^{}_{\mrm{poly}}(\K[[\bsym{t}]])$ and
$\mcal{D}^{}_{\mrm{poly}}(\K[[\bsym{t}]])$
as topological DG Lie algebras. This can certainly be done: 
put on $\mcal{T}^{p}_{\mrm{poly}}(\K[[\bsym{t}]])$ and
$\mrm{F}_m \mcal{D}^{p}_{\mrm{poly}}(\K[[\bsym{t}]])$
the $\bsym{t}$-adic topology, and put on 
$\mcal{D}^{p}_{\mrm{poly}}(\K[[\bsym{t}]]) = \lim_{m \to}$ \linebreak
$\mrm{F}_m \mcal{D}^{p}_{\mrm{poly}}(\K[[\bsym{t}]])$
the direct limit topology (see \cite[Section 1.1]{Ye1}). 
However the dir-inv structure is better suited for our work.
\end{rem}

\begin{dfn}
For $p \geq 0$ let $\mcal{D}^{\mrm{nor}, p}_{\mrm{poly}}(C)$ 
be the submodule of 
$\mcal{D}^p_{\mrm{poly}}(C)$ consisting of poly differen\-tial 
operators $\phi$ such that $\phi(c_1, \ldots, c_{p+1}) = 0$
if $c_i = 1$ for some $i$. For $p = -1$ we let
$\mcal{D}^{\mrm{nor}, -1}_{\mrm{poly}}(C) := C$. 
Define
$\mcal{D}^{\mrm{nor}}_{\mrm{poly}}(C) := 
\boplus_{p \geq -1} \, \mcal{D}^{\mrm{nor}, p}_{\mrm{poly}}(C)$. 
We call $\mcal{D}^{\mrm{nor}}_{\mrm{poly}}(C)$ the algebra of {\em 
normalized poly differential operators}.
\end{dfn}

 From the formulas for the Gerstenhaber bracket and the 
Hochschild differential (see \cite[Section 3.4.2]{Ko1})
it immediately follows that 
$\mcal{D}^{\mrm{nor}}_{\mrm{poly}}(C)$ is a sub DG Lie algebra of
$\mcal{D}_{\mrm{poly}}(C)$.

For any integer $p \geq 1$ there is a $C$-linear homomorphism
\[ \mcal{U}_1 : 
\mcal{T}^{p-1}_{\mrm{poly}}(C) \to 
\mcal{D}^{\mrm{nor}, p-1}_{\mrm{poly}}(C) \]
with formula 
\begin{equation} \label{eqn6.1}
 \mcal{U}_1(\xi_1 \wedge \cdots \wedge \xi_p)
(c_1, \ldots, c_p) := 
{\smfrac{1}{p!}} \sum_{\sigma \in \mfrak{S}_p} \opn{sgn}(\sigma)
\xi_{\sigma(1)}(c_1) \cdots \xi_{\sigma(p)}(c_p) 
\end{equation}
for elements $\xi_1, \ldots, \xi_p \in \mcal{T}_C$
and
$c_1, \ldots, c_p \in C$. For $p = 0$ the map
$\mcal{U}_1 : \mcal{T}^{-1}_{\mrm{poly}}(C) \to 
\mcal{D}^{\mrm{nor}, -1}_{\mrm{poly}}(C)$
is the identity (of $C$).

Suppose $M$ and $N$ are complexes in 
$\cat{Dir} \cat{Inv} \cat{Mod} C$ and $\phi, \phi' : M \to N$ are 
morphisms of complexes in $\cat{Dir} \cat{Inv} \cat{Mod} C$
(i.e.\ all maps are continuous for the dir-inv structures).
We say $\phi$ and $\phi'$ are
homotopic if there is a degree $-1$ homomorphism of graded
dir-inv modules $\eta : M \to N$ such that
$\d_N \circ \eta + \eta \circ \d_M = \phi - \phi'$.
We say that $\phi : M \to N$ is a homotopy equivalence in 
$\cat{Dir} \cat{Inv} \cat{Mod} C$ if there is a morphism of 
complexes $\psi : N \to M$ in $\cat{Dir} \cat{Inv} \cat{Mod} C$
such that
$\psi \circ \phi$ is homotopic to $\bsym{1}_M$ and
$\phi \circ \psi$ is homotopic to $\bsym{1}_N$.

\begin{thm} \label{thm6.1} 
Let $C$ be a commutative $\K$-algebra with ideal $\mfrak{c}$. 
Assume $C$ is noetherian and $\mfrak{c}$-adically complete. 
Also assume there is a $\K$-algebra homomorphism 
$\K[t_1, \ldots, t_n] \to C$ which is flat and 
$\mfrak{c}$-adically formally \'etale. Then the homomorphism
$\mcal{U}_1 : \mcal{T}_{\mrm{poly}}(C) \to 
\mcal{D}^{\mrm{nor}}_{\mrm{poly}}(C)$
and the inclusion 
$\mcal{D}^{\mrm{nor}}_{\mrm{poly}}(C) \to \mcal{D}_{\mrm{poly}}(C)$
are both homotopy equivalences in
$\cat{Dir} \cat{Inv} \cat{Mod} C$.
\end{thm}

\begin{proof}
Recall that 
$\mcal{B}_{q}(C) = \mcal{B}^{-q}(C) := C^{\otimes (q+2)}$,
and this is a $\mcal{B}_{0}(C)$-algebra via the extreme factors. 
So
$\mcal{B}_{q}(C) \cong \mcal{B}_{0}(C) \otimes C^{\otimes q}$
as $\mcal{B}_{0}(C)$-modules. 
Let $\ol{C} := C / \K$, the quotient $\K$-module, and define
$\mcal{B}_{q}^{\mrm{nor}}(C) = \mcal{B}^{\mrm{nor}, -q}(C) 
:= \mcal{B}_{0}(C) \otimes \ol{C}^{\otimes q}$,
the $q$-th normalized bar module of $C$. 
According to \cite[Section X.2]{ML}, 
$\mcal{B}^{\mrm{nor}}(C) :=
\boplus_q \, \mcal{B}^{\mrm{nor}, -q}(C)$
has a coboundary operator such that the obvious surjection
$\phi : \mcal{B}(C) \to \mcal{B}^{\mrm{nor}}(C)$
is a quasi-isomorphism of complexes of $\mcal{B}^0(C)$-modules.

Define
\[ \mcal{C}_{q}^{\mrm{nor}}(C) = \mcal{C}^{\mrm{nor}, -q}(C) 
:= C \otimes_{\mcal{B}_{0}(C)} \mcal{B}^{\mrm{nor}}_{q}(C) \cong
C \otimes \ol{C}^{\otimes q} . \]
Because the complexes $\mcal{B}(C)$ and $\mcal{B}^{\mrm{nor}}(C)$
are bounded above and consist of free $\mcal{B}_{0}(C)$-modules, 
it follows that 
$\phi : \mcal{C}_{}(C) \to \mcal{C}_{}^{\mrm{nor}}(C)$
is a quasi-isomorphism of complexes of $C$-modules.
Let $\what{\Omega}^q_C$ be the $\mfrak{c}$-adic completion of 
$\Omega^q_C$, so that 
$\what{\Omega}^q_C \cong C \otimes_{\K[\bsym{t}]} 
\Omega^q_{\K[\bsym{t}]}$.
There is a $C$-linear homomorphism
$\psi : \mcal{C}_{q}^{\mrm{nor}}(C) \to \Omega^q_C$
with formula
\[ \psi(1 \otimes (c_1 \otimes \cdots \otimes c_q)) :=
\d(c_1) \wedge \cdots \wedge \d(c_q) . \]

Consider the polynomial algebra 
$\K[\bsym{t}] = \K[t_1, \ldots, t_n]$.
For $i \in \{1, \ldots, n\}$ and $j \in \{1, \ldots, q\}$ let
\[ \til{\d}_j (t_i) := 
\underset{j}{\underbrace{1 \otimes \cdots \otimes 1}}
\otimes (t_i \otimes 1 - 1 \otimes t_i) \otimes 1 
\otimes \cdots \otimes 1 \in \mcal{B}_{q}(\K[\bsym{t}]) , \] 
and use the same expression to denote the image of this element in 
$\mcal{C}_{q}(\K[\bsym{t}])$.
It is easy to verify that $\mcal{C}_{q}(\K[\bsym{t}])$ is a 
polynomial algebra over $\K[\bsym{t}]$ in the set of generators 
$\{ \til{\d}_j (t_i) \}$.
Another easy calculation shows that 
$\opn{Ker} \big( \phi : \mcal{C}_{q}(\K[\bsym{t}]) \to
\mcal{C}^{\mrm{nor}}_{q}(\K[\bsym{t}]) \big)$
is generated as $\K[\bsym{t}]$-module by monomials 
in elements of the set $\{ \til{\d}_j (t_i) \}$.

Let's introduce a grading on $\mcal{C}_{q}(\K[\bsym{t}])$ by
$\opn{deg}(\til{\d}_j (t_i)) := 1$ and $\opn{deg}(t_i) := 0$.
The coboundary operator of $\mcal{C}_{}(\K[\bsym{t}])$ 
has degree $0$ in this grading.
The grading is inherited by 
$\mcal{C}^{\mrm{nor}}_{q}(\K[\bsym{t}])$, and hence
$\phi : \mcal{C}(\K[\bsym{t}]) \to
\mcal{C}^{\mrm{nor}}_{}(\K[\bsym{t}])$
is a quasi-isomorphism of complexes in $\cat{GrMod} \K[\bsym{t}]$,
the category of graded $\K[\bsym{t}]$-modules. 
Also let's put a grading on $\Omega^q_{\K[\bsym{t}]}$
with $\opn{deg}(\d(t_i)) := 1$. By \cite[Lemma 4.3]{Ye2}, 
$\psi : \mcal{C}^{\mrm{nor}}(\K[\bsym{t}]) \to
\boplus_q \, \Omega^q_{\K[\bsym{t}]}[q]$
is a quasi-isomorphism in $\cat{GrMod} \K[\bsym{t}]$.
Because we are dealing with bounded above complexes of free graded 
$\K[\bsym{t}]$-modules it follows that both $\phi$ and $\psi$ are 
homotopy equivalences in $\cat{GrMod} \K[\bsym{t}]$.

Now let's go back to the formally \'etale homomorphism
$\K[\bsym{t}] \to C$. We get homotopy equivalences
\[ C \otimes_{\K[\bsym{t}]} \mcal{C}(\K[\bsym{t}])
\xar{\phi} 
C \otimes_{\K[\bsym{t}]} \mcal{C}^{\mrm{nor}}(\K[\bsym{t}])
\xar{\psi} \boplus_q \, \what{\Omega}^q_{C}[q] \]
in $\cat{GrMod} C$. 
We know that $\what{\mcal{C}}_q(C)$ is a power series algebra in the 
set of generators $\{ \til{\d}_j (t_i) \}$; see \cite[Lemma 2.6]{Ye2}.
Therefore $\what{\mcal{C}}_q(C)$ is isomorphic to the completion 
of
$C \otimes_{\K[\bsym{t}]} \mcal{C}_q(\K[\bsym{t}])$
with respect to the grading (see Example \ref{exa1.4}). 
Define
$\what{\mcal{C}}^{\mrm{nor}}_q(C)$ to be the completion of 
$C \otimes_{\K[\bsym{t}]} \mcal{C}^{\mrm{nor}}_q(\K[\bsym{t}])$
with respect to the grading.
We then have a homotopy equivalence of complexes
in $\cat{Inv} \cat{Mod} C$
\[ \what{\mcal{C}}(C) \to \what{\mcal{C}}^{\mrm{nor}}(C) \to
\boplus_q \, \what{\Omega}^q_{C}[q] . \]
Applying $\opn{Hom}^{\mrm{cont}}_C(-, C)$ we arrive at 
quasi-isomorphisms
\[ \boplus_q \, (\bwedge^q_{C} \mcal{T}_C)[-q] \to
\mcal{C}^{\mrm{nor}}_{\mrm{cd}}(C) \to
\mcal{C}^{}_{\mrm{cd}}(C) , \]
where by definition
$\mcal{C}^{\mrm{nor}}_{\mrm{cd}}(C)$ is the continuous dual of 
$\what{\mcal{C}}^{\mrm{nor}}(C)$.
An easy calculation shows that
$\mcal{C}^{\mrm{nor}, q}_{\mrm{cd}}(C) = 
\mcal{D}^{\mrm{nor}, q-1}_{\mrm{poly}}(C)$.
\end{proof}

One instance to which this theorem applies is 
$C := \K[[t_1, \ldots, t_n]]$. Here is another:

\begin{cor}
Suppose $C$ is a smooth $\K$-algebra. Then  the homomorphism
$\mcal{U}_1 : \mcal{T}_{\mrm{poly}}(C) \to 
\mcal{D}^{\mrm{nor}}_{\mrm{poly}}(C)$
and the inclusion 
$\mcal{D}^{\mrm{nor}}_{\mrm{poly}}(C) \to \mcal{D}_{\mrm{poly}}(C)$
are both quasi-isomorphisms.
\end{cor}

\begin{proof}
There is an open covering
$\opn{Spec} C = \bigcup \opn{Spec} C_i$
such that for every $i$ there is an \'etale homomorphism
$\K[t_1, \ldots, t_n] \to C_i$. Now use Theorem \ref{thm6.1}.
\end{proof}

Here is a slight variation of the celebrated result 
of Kontsevich, known as the {\em Formality Theorem}
\cite[Theorem 6.4]{Ko1}.

\begin{thm} \label{thm2.0}
Let $\K[\bsym{t}] = \K[t_1, \ldots, t_n]$ 
be the polynomial algebra in $n$ variables, and 
assume that $\mbb{R} \subset \K$. There is a collection of 
$\K$-linear homomorphisms 
\[ \mcal{U}_j :\, \bwedge^{j}
\mcal{T}_{\mrm{poly}}(\K[\bsym{t}]) \to
\mcal{D}_{\mrm{poly}}(\K[\bsym{t}]) , \]
indexed by $j \in \{ 1, 2, \ldots \}$, satisfying the 
following conditions.
\begin{enumerate}
\rmitem{i} The sequence $\mcal{U} = \{ \mcal{U}_j \}$ 
is an $\mrm{L}_{\infty}$-morphism
$\mcal{T}_{\mrm{poly}}(\K[\bsym{t}]) \to
\mcal{D}_{\mrm{poly}}(\K[\bsym{t}])$.
\rmitem{ii} Each $\mcal{U}_j$ is a poly differential operator 
of $\K[\bsym{t}]$-modules.
\rmitem{iii} Each $\mcal{U}_j$ is equivariant for the standard 
action of $\mrm{GL}_n(\K)$ on $\K[\bsym{t}]$.
\rmitem{iv} The homomorphism $\mcal{U}_1$ is given by equation 
\tup{(\ref{eqn6.1})}.
\rmitem{v} For any $j \geq 2$ and
$\alpha_1, \ldots, \alpha_j \in 
\mcal{T}^{0}_{\mrm{poly}}(\K[\bsym{t}])$
one has
$\mcal{U}_j(\alpha_1 \wedge \cdots \wedge \alpha_j) = 0$.
\rmitem{vi} For any $j \geq 2$,
$\alpha_1 \in \mfrak{gl}_n(\K) \subset 
\mcal{T}^{0}_{\mrm{poly}}(\K[\bsym{t}])$
and
$\alpha_2, \ldots, \alpha_j \in 
\mcal{T}^{}_{\mrm{poly}}(\K[\bsym{t}])$
one has
$\mcal{U}_j(\alpha_1 \wedge \cdots \wedge \alpha_j) = 0$.
\end{enumerate}
\end{thm}

\begin{proof}
First let's assume that $\K = \mbb{R}$.
Theorem 6.4 in \cite{Ko1} talks about the differentiable manifold 
$\mbb{R}^n$, and considers $\mrm{C}^{\infty}$ functions on it, 
rather than polynomial functions. However, by construction
the operators $\mcal{U}_j$ are multi differential operators
with polynomial coefficients (see \cite[Section 6.3]{Ko1}).
Therefore they descend to operators
\[ \mcal{U}_j : \bwedge^{j} 
\mcal{T}^{}_{\mrm{poly}}(\mbb{R}[\bsym{t}]) \to
\mcal{D}^{}_{\mrm{poly}}(\mbb{R}[\bsym{t}]) , \]
and conditions (i) and (ii) hold. Conditions (iii), (v)
and (vi) are properties P3, P4 and P5 
respectively in \cite[Section 7]{Ko1}. For
condition (iv) see \cite[Sections 4.6.1-2]{Ko1}.

For a field extension $\mbb{R} \subset \K$ use base change. 
\end{proof}

\begin{rem} \label{rem8.1}
It is likely that the operator $\mcal{U}_j$ sends
$\bwedge^{j} \mcal{T}_{\mrm{poly}}(\K[\bsym{t}])$
into \linebreak $\mcal{D}^{\mrm{nor}}_{\mrm{poly}}(\K[\bsym{t}])$. 
This is clear for $j = 1$, where 
$\mcal{U}_1(\mcal{T}_{\mrm{poly}}(\K[\bsym{t}])) = 
\mrm{F}_1 \mcal{D}^{\mrm{nor}}_{\mrm{poly}}(\K[\bsym{t}])$;
but this requires checking for $ j \geq 2$.
\end{rem}

In the next theorem $\mcal{T}_{\mrm{poly}}(\K[[\bsym{t}]])$ and
$\mcal{D}_{\mrm{poly}}(\K[[\bsym{t}]])$ are considered as DG Lie 
algebras in 
$\cat{Dir} \cat{Inv} \cat{Mod} \K$, 
with their $\bsym{t}$-adic dir-inv structures. 
Recall the notions of twisted DG Lie algebra (Lemma \ref{lem9.6})
and multilinear extensions of $\mrm{L}_{\infty}$ morphisms
(Proposition \ref{prop10.1}).

\begin{thm} \label{thm4.10}
Assume $\mbb{R} \subset \K$.
Let $A = \boplus_{i \geq 0} A^i$ be a complete super-commutative 
associative unital DG algebra in 
$\cat{Dir} \cat{Inv} \cat{Mod} \K$. Consider the 
induced continuous $A$-multilinear $\mrm{L}_{\infty}$ morphism
\[ \mcal{U}_A : 
A \hatotimes{} \mcal{T}_{\mrm{poly}}(\K[[\bsym{t}]]) \to
A \hatotimes{} \mcal{D}_{\mrm{poly}}(\K[[\bsym{t}]]) . \]
Suppose 
$\omega \in 
A^1 \hatotimes{} \mcal{T}^0_{\mrm{poly}}(\K[[\bsym{t}]])$
is a solution of the Maurer-Cartan equation in 
$A \hatotimes{} \mcal{T}_{\mrm{poly}}(\K[[\bsym{t}]])$.
Define 
$\omega' := (\partial^1 \mcal{U}_A)(\omega) \in
A^1 \hatotimes{} \mcal{D}^0_{\mrm{poly}}(\K[[\bsym{t}]])$.
Then $\omega'$ is a solution of the the Maurer-Cartan equation in 
$A \hatotimes{} \mcal{D}_{\mrm{poly}}(\K[[\bsym{t}]])$,
and there is continuous $A$-multilinear $\mrm{L}_{\infty}$ 
quasi-isomorphism 
\[ \mcal{U}_{A, \omega} :
\big( A \hatotimes{} \mcal{T}_{\mrm{poly}}(\K[[\bsym{t}]]) 
\big)_{\omega} \to
\big( A \hatotimes{} \mcal{D}_{\mrm{poly}}(\K[[\bsym{t}]]) 
\big)_{\omega'} \]
whose Taylor coefficients are
\[ (\partial^j \mcal{U}_{A, \omega})(\alpha) :=
\sum_{k \geq 0} \smfrac{1}{(j+k)!} 
(\partial^{j+k} \mcal{U}_{A})(\omega^k \wedge \alpha)  \]
for 
$\alpha \in \sprod^j
\bigl( A \hatotimes{} \mcal{T}_{\mrm{poly}}(\K[[\bsym{t}]])
\bigr)$.
\end{thm}

\begin{proof}
By condition (ii) of Theorem \ref{thm2.0}, and by Proposition 
\ref{prop1.9}, each operator
$\partial^j \mcal{U} := \mcal{U}_{j}$ 
is continuous for the $\bsym{t}$-adic dir-inv structures on 
$\mcal{T}_{\mrm{poly}}(\K[[\bsym{t}]])$ and
$\mcal{D}_{\mrm{poly}}(\K[[\bsym{t}]])$.
Therefore there is a unique continuous $A$-multilinear extension
$\partial^j \mcal{U}_A$. 
Condition (v) of Theorem \ref{thm2.0} implies that 
$\partial^j \mcal{U}_A(\omega^j) = 0$ for $j \geq 2$. 
By Theorem \ref{thm9.2} we get an 
$\mrm{L}_{\infty}$ morphism $\mcal{U}_{A, \omega}$.

It remains to prove that $\partial^1 \mcal{U}_A$ is a 
quasi-isomorphism. According to Theorem \ref{thm6.1} 
For every $i$ the $\K$-linear homomorphism
\[ \partial^1 \mcal{U}_A :
A^i \hatotimes{} \mcal{T}_{\mrm{poly}}(\K[[\bsym{t}]])
\to 
A^i \hatotimes{} \mcal{T}_{\mrm{poly}}(\K[[\bsym{t}]]) \]
is a quasi-isomorphism. Since we are looking at bounded below 
complexes, a spectral sequence argument implies that 
\[ \partial^1 \mcal{U}_A :
A \hatotimes{} \mcal{T}_{\mrm{poly}}(\K[[\bsym{t}]])
\to 
A \hatotimes{} \mcal{T}_{\mrm{poly}}(\K[[\bsym{t}]]) \]
is a quasi-isomorphism.
\end{proof}



\begin{thebibliography}{EGA ~IV}
\bibitem[CFT]{CFT} S. Cattaneo, G. Felder and L. Tomassini, 
    From local to global deformation quantization of Poisson 
    manifolds, Duke Math.\ J.\ \textbf{115} (2002), 
    no.\ 2, 329-352.
\bibitem[EGA ~I]{EGA-I} A.\ Grothendieck and J.\ Dieudonn\'{e},
     ``\'{E}l\'{e}ments de G\'{e}ometrie Alg\'{e}brique I,''
     Springer, Berlin, 1971.
\bibitem[EGA ~IV]{EGA-IV} A.\ Grothendieck and J.\ Dieudonn\'{e}, 
    ``\'{E}l\'{e}ments de G\'{e}ometrie Alg\'{e}brique IV,''
    Publ.\ Math.\ IHES \textbf{32} (1967).
\bibitem[Fu]{Fu} K. Fukaya, 
    Deformation theory, homological     algebra, and mirror 
    symmetry, in
    ``Geometry and Physics of Branes'' (Como, 2001), pp.\ 121-209, 
    Ser.\ High Energy Phys.\ Cosmol.\ Gravit. (IOP, Bristol, 2003).  
\bibitem[Ge]{Ge} M. Gerstenhaber,
    On the deformation of rings and algebras,
    Ann.\ of Math.\ {\bf 79} (1964), 59-103.
\bibitem[GK]{GK} I.M. Gelfand and D.A. Kazhdan, 
    Some problems of differential geometry and the calculation of 
    cohomologies of Lie algebras of vector fields, 
    Soviet Math.\ Dokl.\ \textbf{12} (1971), no.\ 5, 1367-1370. 
\bibitem[Ke]{Ke} B. Keller, 
    Introduction to Kontsevich's quantization theorem,
    preprint.
\bibitem[Ko1]{Ko1} M. Kontsevich,
    Deformation quantization of Poisson manifolds,
    Lett.\ Math.\ Phys.\ {\bf 66} (2003), no.\ 3, 157-216.
\bibitem[Ko2]{Ko2} M. Kontsevich,
    Operads and Motives in deformation quantization, 
    Lett.\ Math.\ Phys.\ {\bf 48} (1999), 35-72.
\bibitem[Ko3]{Ko3} M. Kontsevich,
    Deformation quantization of algebraic varieties,
    Lett.\ Math.\ Phys.\ {\bf 56} (2001), no.\ 3, 271-294.
\bibitem[Le]{Le} K. Lef\'evre-Hasegawa,
    ``Sur les $\mrm{A}_{\infty}$-Cat\'egories,'' thesis.
\bibitem[ML]{ML} S. MacLane, ``Homology,''
    Reprint of the 1975 edition, Springer-Verlag.
\bibitem[ScSt]{ScSt} M. Schlessinger and J. Stasheff,
   The Lie algebra structure of tangent cohomology and deformation
   theory, J. Pure Appl.\ Algebra {\bf 38} (1985), no. 2-3, 313-322.
\bibitem[Ye1]{Ye1} A.\ Yekutieli,    
    The Continuous Hochschild Cochain Complex of a Scheme,
    Canadian J. Math.\ \textbf{54} (2002), 1319-1337. 
\bibitem[Ye2]{Ye2} A.\ Yekutieli,    
    Deformation Quantization in Algebraic Geometry,
    eprint math.AG/0310399  at http://arxiv.org. 
\end{thebibliography}
\end{document}